\documentclass[11pt]{article}

\usepackage{amsmath,amssymb,amsthm,mathtools}
\usepackage[english]{babel}

\usepackage[letterpaper,top=2cm,bottom=2cm,left=3cm,right=2.8cm,marginparwidth=1.75cm]{geometry}

\usepackage{booktabs}

\usepackage{subcaption}
\usepackage{amsmath}
\usepackage{graphicx}
\usepackage[colorlinks=true, allcolors=blue]{hyperref}
\usepackage[ruled]{algorithm2e}
\usepackage{enumitem}
\numberwithin{equation}{section}
\newtheorem{theorem}{Theorem}[section]
\newtheorem{lemma}[theorem]{Lemma}

\newtheorem{assumption}[theorem]{Assumption}

\newtheorem{definition}[theorem]{Definition}

\newcommand\keywords[1]{\textbf{Keywords}: #1}

\title{NeRI: Neural Residual Iteration for Elliptic Operator Equations\\ and Variational Inequalities}

\author{Michael Hinterm\"uller$^{1,2}$ \quad Jianfeng Ning$^{1}$\\[0.5em]
\small $^{1}$Weierstrass Institute for Applied Analysis and Stochastics, Berlin, Germany\\
\small $^{2}$Institute for Mathematics, Humboldt-Universit\"at zu Berlin, Berlin, Germany\\[0.4em]
\small \texttt{hintermueller@wias-berlin.de} \quad \texttt{ning@wias-berlin.de}}

\date{}
\begin{document}
\maketitle

\begin{abstract}
In this work, we propose a so-called {\it Neural Residual Iteration} (NeRI), a weak-form neural-network framework for solving a class of elliptic operator equations and variational inequalities, including nonlinear and nonsymmetric problems that may not admit a global energy functional. The key idea is to construct, at each iteration, a state-dependent linear functional whose derivative with respect to the trial function recovers the current weak residual. A single gradient step is then performed on the residual, yielding a residual-driven neural update without requiring a global energy functional or an additional test network. From a function-space perspective, NeRI can be interpreted as a neural-tangent-kernel-preconditioned realization of classical residual iterations, thereby bridging classical iterative solvers and neural-network optimization. We further extend NeRI to elliptic variational inequalities (VIs) through penalized and constraint-preserving formulations. For wide two-layer networks under gradient flow, we prove finite-time lazy-training estimates and derive an $H^1$-error decomposition relative to the exact weak solution. Numerical experiments on several elliptic problems demonstrate the effectiveness of NeRI, particularly for nonsymmetric and low-regularity problems.

\end{abstract}

\keywords{
neural residual iteration, elliptic operator equations, variational inequalities,
neural tangent kernel, lazy training}

\section{Introduction}

Partial differential equations (PDEs), and elliptic equations in particular, play a fundamental role in modeling a wide range of stationary phenomena in science and engineering. Their accurate and efficient numerical solution has therefore been a central problem in scientific computing. Classical numerical methods, including finite difference, finite element, finite volume, and spectral methods, have been extensively studied for elliptic equations and have achieved great success in both theory and applications; see, e.g., \cite{ames2014numerical,dauge2006elliptic} and the references therein. In particular, variational formulations and finite element methods provide a powerful framework for general elliptic boundary value problems and extend naturally to elliptic variational inequalities arising from unilateral or convex constraints, such as obstacle and contact problems \cite{kinderlehrer2000introduction,tremolieres2011numerical}. However, traditional mesh-based methods may become computationally expensive for high-dimensional problems, complex geometries, and multiscale structures.

Neural-network-based methods offer an alternative by representing the unknown solution as a trainable function rather than constructing it on a prescribed mesh. Physics-informed neural networks (PINNs) \cite{raissi2019physics} incorporate the differential equation, boundary conditions, and available data into a residual-based training objective, using automatic differentiation \cite{baydin2018automatic} to evaluate the governing equations. These strong-form methods are flexible, but may require expensive evaluation of higher-order derivatives \cite{lu2021deepxde} and can be less effective for solutions with limited regularity \cite{hu2024solving,hu2023solving}.

For elliptic problems, weak and variational formulations provide an important alternative. Compared with strong-form residual methods, weak-form approaches usually require only lower-order derivatives of the neural trial function and are therefore better suited to solutions with limited regularity. The deep Ritz method \cite{yu2018deep} uses an energy functional associated with the PDE as the loss function and applies to elliptic problems admitting a variational potential, including linear symmetric problems and certain nonlinear potential equations. Its applicability is limited, however, when the underlying operator is non-potential and no suitable global energy functional is available. To go beyond this setting, the iterative deep Ritz method (IDRM) \cite{hu2025iterative} treats monotone elliptic problems, including non-potential operators, through a sequence of surrogate variational problems. Another weak-form approach is the weak adversarial network (WAN) \cite{zang2020weak}, which parameterizes both trial and test functions and casts the PDE into an operator-norm or minimax formulation. These approaches broaden the applicability of neural solvers, but also introduce additional computational or optimization challenges: IDRM requires solving a sequence of surrogate neural minimization problems, while WAN typically involves two networks and a minimax optimization. Related developments and numerical analyses of neural PDE methods can be found, for example, in \cite{uriarte2023deep,gao2023failure,lee2026two,weng2026deep,jin2026convergence,jiao2025drm,de2024numerical,wang2022and}. Data-driven operator-learning approaches, such as DeepONet \cite{lu2021learning,sojitra2026fedonet} and Fourier neural operators \cite{lifourier,li2025d}, provide another related direction; see \cite{nganyu2023deep,toscano2025pinns,luo2025physics} for recent reviews.

Recently, neural methods have also been developed for elliptic variational inequalities, where constraints and potentially nonsmooth solutions introduce further challenges. Existing approaches include energy-based and penalty formulations for obstacle problems \cite{zhao2022two,cheng2023deep}, weak adversarial methods based on regularized gap functions \cite{alphonse2024neural}, and Prox-PINNs based on proximal reformulations \cite{gao2025prox}. In particular, the weak adversarial approach in \cite{alphonse2024neural} accommodates nonsymmetric operators, while Prox-PINNs \cite{gao2025prox} use proximal operators to reformulate VIs as nonlinear equations. 

Alongside these methodological developments, the neural tangent kernel (NTK) provides a framework for understanding neural-network training through function-space dynamics \cite{jacot2018neural}. Lazy-training theory relates these dynamics to those of the model linearized at initialization \cite{chizat2019lazy}. For PINNs, Wang et al.~\cite{wang2022and} use the NTK to analyze discrepancies in the convergence rates of different loss components. Jin and Wu \cite{jin2026convergence} establish empirical-loss convergence for wide two-layer PINNs trained by stochastic gradient methods for the Poisson equation, relying on the positive definiteness of suitable Gram matrices. Related kernel-based work studies Sobolev-norm accuracy and statistical optimality for learning elliptic equations \cite{lu2022sobolev}.

In this work, we propose \emph{Neural Residual Iteration} (NeRI), a framework based directly on the weak operator residual. For an elliptic operator equation $\mathcal{A}(u)=f$, we construct, at each iterate $u_k$, a state-dependent linear functional whose derivative with respect to the trial function recovers $\mathcal{A}(u_k)-f$. In the neural parameterization, a single gradient step is then performed, and the residual is re-evaluated at the updated approximation. Unlike IDRM, the state-dependent functional in NeRI is not treated as a surrogate variational problem to be minimized; it is used only to generate the current residual-driven update. Thus, NeRI has no inner minimization loop and is formulated directly as a neural realization of residual iteration. In the potential case, its update direction coincides with the deep Ritz gradient direction, while the construction remains applicable to nonlinear and nonsymmetric non-potential operators. We also provide a weak adversarial interpretation through an implicitly defined test function, without introducing an additional network. From a function-space perspective, NeRI induces, to first order, a neural-tangent-kernel-preconditioned residual iteration, connecting classical operator iterations with neural-network optimization \cite{ortega2000iterative,zeidler2013nonlinear}. A closely related non-gradient method is proposed by Peng et al.~\cite{peng2023non}, which updates network parameters using strong-form PDE residuals evaluated by finite differences. NeRI shares this residual-driven principle but uses weak residual pairings involving only first-order spatial derivatives for second-order divergence-form problems. This formulation naturally accommodates low-regularity weak solutions without requiring pointwise differential residuals. Our work further extends the framework to elliptic variational inequalities and develops a Sobolev-space analysis.

An important feature of our analysis is that the neural tangent kernel is treated as an operator from the dual of the Sobolev trial space into the trial space itself. Under a gradient-flow formulation, we prove finite-time lazy-training estimates and a frozen-kernel approximation. For wide two-layer networks, the kernel drift is $O_{\mathbb{P}}((\log m)/\sqrt{m})$ under suitable activation and initialization assumptions. In contrast to the empirical-loss convergence analysis in \cite{jin2026convergence}, we combine these estimates to bound, on fixed finite time intervals, the $H^1(\Omega)$-error relative to the exact weak solution, separating finite-width, frozen-flow optimization, and reduced-space approximation errors. 

The remainder of the paper is organized as follows. Section~\ref{sec:problem_setting} introduces the operator setting and classical residual iteration. Section~\ref{sec:NeRI} presents NeRI, its function-space interpretation, and the treatment of Dirichlet boundary conditions. Section~\ref{sec:VI} extends the framework to elliptic VIs. Section~\ref{sec:lazy} develops the lazy-training and error analysis. Numerical experiments and concluding remarks are presented in Sections~\ref{sec:numerical} and~\ref{sec:conclusion}, respectively.

\section{Problem setting and classical residual iteration}
\label{sec:problem_setting}

Let $\Omega\subset\mathbb R^d$ be a bounded domain with boundary $\partial\Omega$, and let $V$ be the Hilbert space used in the weak formulation, with dual $V^*$ and duality pairing $\langle\cdot,\cdot\rangle_{V^*,V}$. For second-order problems, typical choices are $V=H^1(\Omega)$ for natural boundary conditions and $V=H_0^1(\Omega)$ for homogeneous Dirichlet conditions \cite{FournierAdams}. Nonhomogeneous Dirichlet data can be handled by a lifting $u=\hat g_D+w$, where $\hat g_D\in H^1(\Omega)$ extends the boundary data and $w\in H_0^1(\Omega)$.

We consider the abstract elliptic operator equation
\begin{equation}\label{equa:ellPDEs}
    \text{find } u\in V:\qquad
    \mathcal A(u)=f
    \quad \text{in } V^*,
\end{equation}
where $\mathcal A:V\to V^*$ is a possibly nonlinear and nonsymmetric operator, and $f\in V^*$ is given. A typical second-order weak form is
\begin{align}
    \langle\mathcal A(u),v\rangle_{V^*,V}
    &=
    \int_\Omega
    \bigl(a_0(x,u,\nabla u)v
    +a_1(x,u,\nabla u)\cdot\nabla v\bigr)\,\mathrm dx
    +\int_{\partial\Omega}b(x,u)v\,\mathrm ds,
    \\
    \langle f,v\rangle_{V^*,V}
    &=
    \int_\Omega
    \bigl(f_0v+f_1\cdot\nabla v\bigr)\,\mathrm dx
    +\int_{\partial\Omega}g_Nv\,\mathrm ds.
\end{align}
The boundary terms account for Neumann or Robin conditions.

Throughout this paper, we assume that $\mathcal A$ is Lipschitz continuous and strongly monotone.

\begin{assumption}\label{assup:Lip+Coer}
There exist constants $L,\gamma>0$ such that, for all $u,v\in V$,
\begin{align}
    \|\mathcal A(u)-\mathcal A(v)\|_{V^*}
    &\le L\|u-v\|_V,
    \\
    \langle\mathcal A(u)-\mathcal A(v),u-v\rangle_{V^*,V}
    &\ge \gamma\|u-v\|_V^2.
    \label{eq:monotone}
\end{align}
\end{assumption}

Under Assumption~\ref{assup:Lip+Coer}, problem~\eqref{equa:ellPDEs} admits a unique solution $u^*$ \cite{zeidler2013nonlinear}. Let $\mathcal R:V^*\to V$ be the Riesz isomorphism defined by $\langle\mathcal Rg,h\rangle_V=\langle g,h\rangle_{V^*,V}$ for $g\in V^*$ and $h\in V$. Given $u_0\in V$ and $\eta>0$, the Richardson-type residual iteration reads
\begin{equation}\label{equa:Richard}
    u_{k+1}
    =u_k-\eta\mathcal R\bigl(\mathcal A(u_k)-f\bigr),
    \qquad k=0,1,2,\ldots.
\end{equation}
The Lipschitz and monotonicity bounds yield
\begin{equation*}
    \|u_{k+1}-u^*\|_V^2
    \le
    \bigl(1-2\eta\gamma+\eta^2L^2\bigr)
    \|u_k-u^*\|_V^2.
\end{equation*}
Hence, the iteration is contractive and converges to $u^*$ whenever $0<\eta<2\gamma/L^2$. More generally, one may consider classical function space iterations of the form $u_{k+1}=u_k-\eta T(u_k),$ where $T(u^*)=0$ and $u\mapsto u-\eta T(u)$ is a contraction for an appropriate choice of the step size $\eta$. This residual-iteration principle is the basis of the neural method developed below.

\section{Neural residual iteration (NeRI)}
\label{sec:NeRI}

\subsection{Motivation and formulation of NeRI}

We first recall the variational case. Suppose that there exists a Fr\'echet-differentiable functional $E:V\to\mathbb R$ such that
\begin{equation}\label{eq:potential_relation}
    D_uE(u)=\mathcal A(u)-f
    \qquad \text{in } V^*.
\end{equation}
Then every solution of $\mathcal A(u)=f$ is a stationary point of $E$. If $E$ is strongly convex, the solution is its unique minimizer. This variational structure forms the basis of the deep Ritz method, which uses $E(u_\theta)$ as the training objective for a neural-network approximation $u_\theta$. For a linear symmetric operator, the energy takes the familiar quadratic form
\begin{equation}\label{eq:ritz_varia}
    E(u)=\frac{1}{2}
    \langle\mathcal A(u),u\rangle_{V^*,V}
    -\langle f,u\rangle_{V^*,V}.
\end{equation}

For general nonlinear or nonsymmetric operators, however, such a global energy need not exist, and a fixed Ritz-type objective may therefore be unavailable. In particular, for a nonsymmetric linear operator, the first variation of \eqref{eq:ritz_varia} involves only the symmetric part of $\mathcal A$ and does not recover the original residual.

The classical residual iteration~\eqref{equa:Richard}, in contrast, does not require $\mathcal A$ to be linear, symmetric, or potential; its contraction relies instead on the Lipschitz continuity and strong monotonicity conditions in Assumption~\ref{assup:Lip+Coer}. Motivated by this observation, NeRI generates a residual-driven neural update without constructing a global energy. In fact, given the current approximation $u_k\in V$ of $u^*$, we define the state-dependent linear functional
\begin{equation}
    J_k(u):=\langle\mathcal A(u_k)-f,u\rangle_{V^*,V}.
\end{equation}
Since $\mathcal A(u_k)-f$ is held fixed when differentiating with respect to $u$, its derivative satisfies $D_uJ_k(u)=\mathcal A(u_k)-f$. Thus, the current residual can be recovered as the derivative of a state-dependent functional even when the operator has no global potential.

The preceding construction is formulated in function space. To obtain a practical neural method, we now restrict the trial functions to a neural-network class and translate the residual direction into a parameter update. For the moment, we consider natural boundary conditions and take $V=H^1(\Omega)$. With a sufficiently regular activation function, the neural realization belongs to $V$, while the natural boundary conditions are incorporated through the weak formulation. Essential Dirichlet conditions are treated separately below in Subsection~\ref{sec:Dirichlet}.

Let $u_\theta\in V$ be a neural-network approximation with parameter vector $\theta\in\mathbb R^p$, and assume that $\theta\mapsto u_\theta$ is Fr\'echet differentiable. At iteration $k$, we define the state-dependent objective
\begin{equation}\label{eq:NeRI_loss}
    \mathcal J_k(\theta):=
    \left\langle
        \operatorname{sg}\bigl(\mathcal A(u_{\theta_k})-f\bigr),
        u_\theta
    \right\rangle_{V^*,V},
\end{equation}
where `$\operatorname{sg}$' denotes the stop-gradient operation: the current residual is treated as a constant when differentiating with respect to $\theta$. In particular, no derivative is propagated through the operator evaluation at $u_{\theta_k}$. We then perform one parameter-gradient step with step size $\eta_k>0$. By the chain rule, this gives
\begin{equation}\label{eq:NeRI_param_update}
        \theta_{k+1}
        =\theta_k-\eta_k
        \nabla_\theta\mathcal J_k(\theta)\big|_{\theta=\theta_k}=\theta_k-\eta_k(D_\theta u_{\theta_k})^*
        \bigl(\mathcal A(u_{\theta_k})-f\bigr),
\end{equation}
where $D_\theta u_{\theta_k}:\mathbb R^p\to V$ is the parameter-to-function derivative and $(D_\theta u_{\theta_k})^*:V^*\to\mathbb R^p$ is its adjoint.

The objective $\mathcal J_k$ is not minimized to convergence; it is used only to generate a single update direction. After the update, the operator is evaluated again at $u_{\theta_{k+1}}$, and a new state-dependent objective is constructed. Unlike a fixed training loss, the objective therefore changes along the iteration and incorporates the current residual information at each step.

For the analysis below, we use plain gradient descent and its continuous-time gradient-flow limit, as is standard in much of the NTK and lazy-training literature~\cite{chizat2019lazy,jacot2018neural,wang2022and,jin2026differential}. In practice, adaptive optimizers such as Adam are widely used for neural PDE solvers because of their favorable empirical optimization behavior; accordingly, we mainly use Adam in Section~\ref{sec:numerical}. The same stop-gradient objective is compatible with Adam, but the resulting history-dependent dynamics are beyond the scope of the present analysis.

This construction also clarifies the connection with the deep Ritz method. Whenever a functional $E$ satisfies \eqref{eq:potential_relation}, the chain rule yields $ \nabla_\theta E(u_\theta)=(D_\theta u_\theta)^*\bigl(\mathcal A(u_\theta)-f\bigr)$, which is exactly the NeRI update direction. Thus, the two methods have the same parameter-gradient direction in the potential case, including the linear symmetric setting \eqref{eq:ritz_varia}. NeRI, however, does not require such an energy and remains applicable to non-potential operators.

\begin{algorithm}[t]
\caption{Neural Residual Iteration (NeRI)}
\label{alg:NeRI}
\KwIn{
A neural network $u_\theta$ with initial parameters $\theta_0$,
step sizes $\{\eta_k\}_{k=0}^{K-1}$, and number of iterations $K$
}
\KwOut{Neural approximation $u_{\theta_K}$}

\For{$k=0,\ldots,K-1$}{
    $\theta_{k+1}=\theta_k-\eta_k
    \nabla_\theta\langle r_k,u_\theta\rangle_{V^*,V}
    \big|_{\theta=\theta_k}$\, \quad where \quad $r_k=\operatorname{sg}\bigl(\mathcal A(u_{\theta_k})-f\bigr)$;
}
\end{algorithm}

Algorithm~\ref{alg:NeRI} summarizes the basic iteration with natural boundary conditions and exact evaluation of the weak residual. Its implementation requires only the weak residual pairing and a standard gradient step, with the residual detached from the computational graph. In practice, the exact duality pairing can be replaced by a quadrature or Monte Carlo approximation.

\subsection{Induced function-space iteration}

We next derive the function-space dynamics induced by NeRI. Although the update is performed in parameter space, it generates a first-order iteration in function space. Assuming that $D_\theta u_\theta$ is locally Lipschitz, a Taylor expansion gives
\begin{equation}
    u_{\theta_{k+1}}
    =
    u_{\theta_k}
    +D_\theta u_{\theta_k}(\theta_{k+1}-\theta_k)
    +\mathcal O(\|\theta_{k+1}-\theta_k\|^2).
\end{equation}
Combining this expansion with \eqref{eq:NeRI_param_update}, we obtain
\begin{equation}
    u_{\theta_{k+1}}
    =
    u_{\theta_k}
    -\eta_k D_\theta u_{\theta_k}
    (D_\theta u_{\theta_k})^*
    \bigl(\mathcal A(u_{\theta_k})-f\bigr)
    +\mathcal O(\eta_k^2).
\end{equation}
The remainders are understood in the $V$-norm. The latter bound is uniform in $k$ provided the Lipschitz constants along the update segments are uniformly bounded and $(D_\theta u_{\theta_k})^*\bigl(\mathcal A(u_{\theta_k})-f\bigr)$ remains bounded.

We define the neural tangent operator at $\theta_k$ by
\begin{equation}
    K_{\theta_k}
    :=
    D_\theta u_{\theta_k}(D_\theta u_{\theta_k})^*
    :V^*\to V.
\end{equation}
Equivalently, writing $\theta=(\theta^{(1)},\ldots,\theta^{(p)})$, we have
\begin{equation*}
    K_{\theta_k}g
    =
    \sum_{i=1}^p
    \langle g,\partial_{\theta^{(i)}}u_{\theta_k}\rangle_{V^*,V}
    \partial_{\theta^{(i)}}u_{\theta_k},
    \qquad g\in V^*.
\end{equation*}
Thus, the induced function-space iteration becomes
\begin{equation}\label{eq:NeRI_function_update}
    u_{\theta_{k+1}}
    =
    u_{\theta_k}
    -\eta_k K_{\theta_k}
    \bigl(\mathcal A(u_{\theta_k})-f\bigr)
    +\mathcal O(\eta_k^2).
\end{equation}
The operator $K_{\theta_k}$ is symmetric in the duality sense and positive semidefinite: for all $g,h\in V^*$,
\begin{equation*}
    \langle g,K_{\theta_k}h\rangle_{V^*,V}
    =
    \langle h,K_{\theta_k}g\rangle_{V^*,V},
    \qquad
    \langle g,K_{\theta_k}g\rangle_{V^*,V}
    =
    |(D_\theta u_{\theta_k})^*g|_{\mathbb R^p}^2
    \ge0.
\end{equation*}
Accordingly, NeRI can be interpreted as a neural-kernel-preconditioned realization of classical residual iterations.

This interpretation is particularly transparent for a linear parameterization $u_\theta=\sum_{i=1}^p\theta^{(i)}\phi_i$, with fixed functions $\phi_i\in V$. In this case, $D_\theta u_\theta$ is independent of $\theta$, and the tangent operator is
\begin{equation}\label{eq:kernel_fix}
    K_\phi g
    =
    \sum_{i=1}^p
    \langle g,\phi_i\rangle_{V^*,V}\phi_i
    =
    \left(\sum_{i=1}^p\phi_i\otimes\phi_i\right)g,
    \qquad g\in V^*.
\end{equation}
The range of the fixed finite rank kernel $K_\phi$ is $V_p:=\operatorname{span}\{\phi_1,\ldots,\phi_p\}$. Since this space is finite-dimensional, $K_\phi$ cannot be uniformly coercive on all of $V^*$ when $V$ is infinite-dimensional. Nevertheless, its positivity on the corresponding reduced space permits a convergence analysis.

To isolate the role of the kernel in residual correction, we first analyze the fixed-kernel iteration. Section~\ref{sec:lazy} then justifies this approximation through finite-time lazy-training estimates for wide two-layer networks. For $u_0\in V$ and $\eta>0$, consider
\begin{equation}\label{equa:ideal_ite}
    u_{k+1}
    =
    u_k-\eta K_\phi\bigl(\mathcal A(u_k)-f\bigr).
\end{equation}
The self-adjoint operator $B:=K_\phi\mathcal R^{-1}:V\to V$ has range $V_p$ and is positive definite on this space. Hence, there exist constants $0<\gamma_K\le M_K$ such that
\begin{equation}\label{eq:K_bounds_fixed}
    \gamma_K\|v\|_V^2
    \le
    \langle v,Bv\rangle_V
    \le
    M_K\|v\|_V^2,
    \qquad v\in V_p.
\end{equation}
Its Moore--Penrose pseudoinverse $B^\dagger$ therefore defines a norm on $V_p$:
\begin{equation*}
    \|v\|_{B^\dagger}^2
    :=
    \langle B^\dagger v,v\rangle_V,
    \qquad v\in V_p.
\end{equation*}
The following result combines a finite-rank preconditioned convergence estimate with a classical C\'ea-type quasi-optimality bound for nonlinear Galerkin approximations; see, e.g., \cite{CongreveWihler2017,GantnerEtAl2018}.

\begin{theorem}\label{thm:fixed_preconditioner}
Suppose Assumption~\ref{assup:Lip+Coer} holds and let $S:=u_0+V_p$. There exists a unique reduced solution $u_S^*\in S$ satisfying
\begin{equation}\label{eq:reduced_problem}
    \langle\mathcal A(u_S^*)-f,v\rangle_{V^*,V}=0,
    \qquad v\in V_p.
\end{equation}
The iterates of \eqref{equa:ideal_ite} satisfy
\begin{equation}\label{eq:finite_rank_descent}
    \|u_{k+1}-u_S^*\|_{B^\dagger}^2
    \le
    \|u_k-u_S^*\|_{B^\dagger}^2
    -\eta\bigl(2\gamma-\eta M_KL^2\bigr)
    \|u_k-u_S^*\|_V^2.
\end{equation}
Consequently, if $0<\eta<2\gamma/(M_KL^2)$, then
\begin{equation}\label{eq:finite_rank_contraction}
    \|u_{k+1}-u_S^*\|_{B^\dagger}^2
    \le
    \rho\|u_k-u_S^*\|_{B^\dagger}^2,\qquad  \rho:=1-\eta\gamma_K(2\gamma-\eta M_KL^2)\in[0,1).
\end{equation}
Moreover, if $u^*$ denotes the solution of $\mathcal A(u^*)=f$, then
\begin{equation}\label{eq:finite_rank_total_error}
    \|u_k-u^*\|_V
    \le
    \sqrt{M_K}\,\rho^{k/2}
    \|u_0-u_S^*\|_{B^\dagger}
    +
    \frac{L}{\gamma}
    \inf_{v\in S}\|v-u^*\|_V.
\end{equation}
\end{theorem}

\begin{proof}
Let $P_p:V\to V_p$ be the orthogonal projection. The reduced operator
$w\mapsto P_p\mathcal R\bigl(\mathcal A(u_0+w)-f\bigr)$
is strongly monotone and Lipschitz on $V_p$, so the reduced problem has a unique solution.

Set $e_k^S:=u_k-u_S^*$ and $r_k^S:=\mathcal A(u_k)-\mathcal A(u_S^*)$. Since $u_0,u_S^*\in S$ and all increments in \eqref{equa:ideal_ite} belong to $V_p$, we have $e_k^S\in V_p$ for all $k\ge0$. By \eqref{eq:reduced_problem},
$K_\phi\bigl(\mathcal A(u_S^*)-f\bigr)=0$, and hence
$e_{k+1}^S=e_k^S-\eta K_\phi r_k^S$.
Using $K_\phi=B\mathcal R$ and $B^\dagger B=P_p$, we obtain
\begin{equation*}
    \|e_{k+1}^S\|_{B^\dagger}^2
    =
    \|e_k^S\|_{B^\dagger}^2
    -2\eta\langle r_k^S,e_k^S\rangle_{V^*,V}
    +\eta^2\langle r_k^S,K_\phi r_k^S\rangle_{V^*,V}.
\end{equation*}
Strong monotonicity gives
$\langle r_k^S,e_k^S\rangle_{V^*,V}\ge\gamma\|e_k^S\|_V^2$,
while the upper spectral bound and Lipschitz continuity imply
\begin{equation*}
    \langle r_k^S,K_\phi r_k^S\rangle_{V^*,V}
    \le
    M_K\|r_k^S\|_{V^*}^2
    \le
    M_KL^2\|e_k^S\|_V^2.
\end{equation*}
This proves \eqref{eq:finite_rank_descent}. Under the stated step-size restriction, substituting
$\|e_k^S\|_V^2\ge\gamma_K\|e_k^S\|_{B^\dagger}^2$
into \eqref{eq:finite_rank_descent} yields
\eqref{eq:finite_rank_contraction}.

We next estimate the reduced-space approximation error. For any $w\in S$, we have $u_S^*-w\in V_p$, so the reduced equation gives
\begin{equation*}
    \langle
        \mathcal A(u_S^*)-\mathcal A(u^*),
        u_S^*-w
    \rangle_{V^*,V}
    =0.
\end{equation*}
Therefore,
\begin{align*}
    \gamma\|u_S^*-u^*\|_V^2
    &\le
    \langle
        \mathcal A(u_S^*)-\mathcal A(u^*),
        u_S^*-u^*
    \rangle_{V^*,V}
    \\
    &=
    \langle
        \mathcal A(u_S^*)-\mathcal A(u^*),
        w-u^*
    \rangle_{V^*,V}
    \\
    &\le
    L\|u_S^*-u^*\|_V\|w-u^*\|_V.
\end{align*}
Taking the infimum over $w\in S$ yields
\begin{equation}\label{eq:finite_rank_cea}
    \|u_S^*-u^*\|_V
    \le
    \frac{L}{\gamma}
    \inf_{v\in S}\|v-u^*\|_V.
\end{equation}
Finally, the upper spectral bound gives
$\|v\|_V^2\le M_K\|v\|_{B^\dagger}^2$ for $v\in V_p$.
Together with \eqref{eq:finite_rank_contraction}, this implies
\begin{equation*}
    \|u_k-u_S^*\|_V
    \le
    \sqrt{M_K}\,\rho^{k/2}
    \|u_0-u_S^*\|_{B^\dagger}.
\end{equation*}
Combining this estimate with \eqref{eq:finite_rank_cea} and the triangle inequality proves \eqref{eq:finite_rank_total_error}.
\end{proof}

\subsection{Treatment of Dirichlet boundary conditions}
\label{sec:Dirichlet}

Unlike natural boundary conditions, essential Dirichlet conditions are not enforced automatically by the weak formulation. We therefore consider two standard treatments within the NeRI framework: the penalty method and the hard-constrained parameterization. We assume that the Dirichlet data $g$ admit a lifting $\hat g\in H^1(\Omega)$ satisfying $\operatorname{Tr}\hat g=g$. On bounded Lipschitz domains, this is equivalent to requiring \(g\in H^{1/2}(\partial\Omega)\).

\subsubsection{Penalty method}

The penalty method augments the NeRI objective with a term that penalizes violations of the prescribed boundary condition
\begin{equation}\label{eq:dirichlet_cond}
    u=g \qquad \text{on } \partial\Omega.
\end{equation}
Let $\operatorname{Tr}:H^1(\Omega)\to L^2(\partial\Omega)$ be the bounded trace operator. The boundary condition is imposed approximately through
\begin{equation*}
    \mathcal P_\beta(u)
    :=
    \frac{\beta}{2}
    \|\operatorname{Tr}u-g\|_{L^2(\partial\Omega)}^2,
\end{equation*}
where $\beta>0$ is the penalty parameter. Its Fr\'echet derivative is $D_u\mathcal P_\beta(u)
=\beta\,\operatorname{Tr}^*\bigl(\operatorname{Tr}u-g\bigr)\in H^1(\Omega)^*$, where $\operatorname{Tr}^*:L^2(\partial\Omega)\to H^1(\Omega)^*$ is the adjoint of the trace operator.

The penalty problem is posed on $V=H^1(\Omega)$, with $\mathcal A:V\to V^*$ and $f\in V^*$ interpreted through the weak formulation on this space. For a neural ansatz $u_\theta=\mathcal N_\theta$, define the penalized NeRI objective at iteration $k$ by
\begin{equation}
    \mathcal J_k^\beta(\theta)
    :=
    \left\langle
        \operatorname{sg}
        \bigl(\mathcal A(u_{\theta_k})-f\bigr),
        u_\theta
    \right\rangle_{V^*,V}
    +\mathcal P_\beta(u_\theta).
\end{equation}
We then perform one parameter-gradient step,
\begin{equation*}
        \theta_{k+1}=
        \theta_k-\eta_k
        \nabla_\theta\mathcal J_k^\beta(\theta_k)
        =
        \theta_k-\eta_k(D_\theta u_{\theta_k})^*
        \Bigl(
            \mathcal A(u_{\theta_k})-f
            +\beta\,\operatorname{Tr}^*
            \bigl(\operatorname{Tr}u_{\theta_k}-g\bigr)
        \Bigr).
\end{equation*}
Introducing the penalized operator
\begin{equation}
    \mathcal A_\beta(u)
    :=
    \mathcal A(u)
    +\beta\,\operatorname{Tr}^*
    \bigl(\operatorname{Tr}u-g\bigr),
\end{equation}
the induced first-order function-space iteration can be written as
\begin{equation*}
    u_{\theta_{k+1}}
    =
    u_{\theta_k}
    -\eta_k K_{\theta_k}
    \bigl(\mathcal A_\beta(u_{\theta_k})-f\bigr)
    +\mathcal O(\eta_k^2),
\end{equation*}
where $K_{\theta_k}=D_\theta u_{\theta_k}(D_\theta u_{\theta_k})^*$. Thus, the penalty method amounts to applying NeRI to $\mathcal A_\beta(u)=f$ while retaining the original neural parameterization. Increasing $\beta$ strengthens the boundary penalty but also increases the Lipschitz bound for $\mathcal A_\beta$, which may reduce the admissible step size and worsen the conditioning of the iteration.

\subsubsection{Hard-constrained parameterization}

A hard-constrained parameterization instead enforces the Dirichlet boundary condition by construction. For \eqref{eq:dirichlet_cond}, one can use an ansatz of the form
\begin{equation}\label{eq:para_boundary}
    u_\theta=\mathcal T(\theta)
    :=
    \hat g+w\mathcal N_\theta,
\end{equation}
where $\hat g\in H^1(\Omega)$ extends the boundary data, with $\operatorname{Tr}\hat g=g$, and $w$ is a boundary factor satisfying
\begin{equation}\label{eq:boundary_func}
    w=0 \quad \text{on } \partial\Omega,
    \qquad
    w>0 \quad \text{in } \Omega.
\end{equation}
Under suitable regularity assumptions, $w\mathcal N_\theta\in H_0^1(\Omega)$, so the prescribed boundary condition holds for every $\theta$. For more general hard-constrained parameterizations, we refer to \cite{hintermuller2026constrained,lu2021physics}.

Using the same ambient space $V=H^1(\Omega)$, the hard-constrained NeRI objective is
\begin{equation*}
    \mathcal J_k^{\mathrm{hc}}(\theta)
    :=
    \left\langle
        \operatorname{sg}
        \bigl(\mathcal A(\mathcal T(\theta_k))-f\bigr),
        \mathcal T(\theta)
    \right\rangle_{V^*,V}.
\end{equation*}
We then perform the gradient update
$\theta_{k+1}=\theta_k-\eta_k
\nabla_\theta\mathcal J_k^{\mathrm{hc}}(\theta)
\big|_{\theta=\theta_k}$.
The induced first-order function-space iteration becomes
\begin{equation*}
    u_{\theta_{k+1}}
    =
    u_{\theta_k}
    -\eta_k\widetilde K_{\theta_k}
    \bigl(\mathcal A(u_{\theta_k})-f\bigr)
    +\mathcal O(\eta_k^2),
\end{equation*}
where
$\widetilde K_{\theta_k}
:=D_\theta\mathcal T(\theta_k)(D_\theta\mathcal T(\theta_k))^*$
is the constrained neural tangent operator. Since
$D_\theta\mathcal T(\theta)=wD_\theta\mathcal N_\theta$
takes values in $H_0^1(\Omega)$, the tangent corrections preserve the boundary condition.

In contrast to the penalty method, the hard-constrained method restricts the trial class to boundary-conforming functions, enforcing the boundary data exactly while changing the effective tangent geometry.

\subsection{Weak adversarial interpretation via an implicit test function}
\label{subsec:implicit_weak_adversarial}

In this subsection, we interpret NeRI from a weak adversarial perspective. We show that its update is gradient-equivalent to a weak residual pairing with an implicitly defined test function associated with the current residual.

Recall that the operator equation can be viewed as minimizing the squared dual residual norm. By the definition of the dual norm,
\begin{equation*}
    \min_{u\in V}\|\mathcal A(u)-f\|_{V^*}^2
    =
    \min_{u\in V}
    \sup_{\varphi\in V,\ \varphi\ne0}
    \frac{
        |\langle\mathcal A(u)-f,\varphi\rangle_{V^*,V}|^2
    }{\|\varphi\|_V^2}.
\end{equation*}
This formulation motivates the weak adversarial network (WAN) \cite{zang2020weak}, which parameterizes both the trial function $u_\theta$ and the test function $\varphi_\eta$. A corresponding objective is
\begin{equation*}
    \mathcal L(u_\theta,\varphi_\eta)
    =
    \frac{
        |\langle\mathcal A(u_\theta)-f,\varphi_\eta\rangle_{V^*,V}|^2
    }{\|\varphi_\eta\|_V^2}.
\end{equation*}
A logarithmic variant, defined when the residual pairing is nonzero, is
\begin{equation*}
    \mathcal L_{\log}(u_\theta,\varphi_\eta)
    =
    \log|\langle\mathcal A(u_\theta)-f,\varphi_\eta\rangle_{V^*,V}|^2
    -\log\|\varphi_\eta\|_V^2.
\end{equation*}
The trial and test networks are then updated alternately through the min--max formulation.

For a nonzero residual, the choice
$\varphi=\mathcal R(\mathcal A(u_\theta)-f)$
attains the supremum and satisfies
\begin{equation*}
    \frac{
        |\langle\mathcal A(u_\theta)-f,\varphi\rangle_{V^*,V}|^2
    }{\|\varphi\|_V^2}
    =
    \|\mathcal A(u_\theta)-f\|_{V^*}^2.
\end{equation*}
Computing this optimal test function generally requires solving a Riesz representation problem. In WAN, it is instead approximated by a separate neural network.

We now relate NeRI to a different implicit test function. Define the unnormalized weak residual pairing
\begin{equation*}
    \widehat{\mathcal L}(u,\varphi)
    :=
    \langle\mathcal A(u)-f,\varphi\rangle_{V^*,V}.
\end{equation*}
Assume that $\mathcal A$ is Fr\'echet differentiable. For fixed $\varphi$, its derivative with respect to $u$ is
$D_u\widehat{\mathcal L}(u,\varphi)
=(D_u\mathcal A(u))^*\varphi$,
where the adjoint is understood in the duality sense. At iteration $k$, the NeRI objective has derivative $\mathcal A(u_k)-f$ with respect to the trial function. We therefore seek a test function $\widehat\varphi_k$ satisfying
\begin{equation*}
    D_u\widehat{\mathcal L}(u_k,\widehat\varphi_k)
    =
    (D_u\mathcal A(u_k))^*\widehat\varphi_k
    =
    \mathcal A(u_k)-f.
\end{equation*}
Under Assumption~\ref{assup:Lip+Coer}, both $D_u\mathcal A(u_k)$ and its adjoint are bounded and coercive, hence isomorphisms from $V$ onto $V^*$. Thus, the implicit test function is well defined by
\begin{equation}\label{eq:implicit_test}
    \widehat\varphi_k
    :=
    \bigl((D_u\mathcal A(u_k))^*\bigr)^{-1}
    \bigl(\mathcal A(u_k)-f\bigr).
\end{equation}
With $\widehat\varphi_k$ held fixed, the derivatives of this pairing and the NeRI objective agree at $u_k$. The corresponding parameter gradients therefore also agree when $u_k=u_{\theta_k}$.

Evaluating the pairing at this implicit test function gives
\begin{equation}\label{eq:implicit_test_residual_norm}
    \widehat{\mathcal L}(u_k,\widehat\varphi_k)
    =
    \left\langle
        \mathcal A(u_k)-f,
        \bigl((D_u\mathcal A(u_k))^*\bigr)^{-1}
        \bigl(\mathcal A(u_k)-f\bigr)
    \right\rangle_{V^*,V}.
\end{equation}
The following lemma shows that this quantity is comparable to the squared dual residual norm.

\begin{lemma}\label{lem:implicit_test_norm_equivalence}
Suppose Assumption~\ref{assup:Lip+Coer} holds and $\mathcal A$ is Fr\'echet differentiable. Then the quantity defined in \eqref{eq:implicit_test_residual_norm} satisfies
\begin{equation}\label{eq:wan_equi_norm}
    \frac{\gamma}{L^2}
    \|\mathcal A(u_k)-f\|_{V^*}^2
    \le
    \widehat{\mathcal L}(u_k,\widehat\varphi_k)
    \le
    \frac{1}{\gamma}
    \|\mathcal A(u_k)-f\|_{V^*}^2.
\end{equation}
\end{lemma}

\begin{proof}
Set $B_k:=D_u\mathcal A(u_k)$. Differentiating the Lipschitz and strong monotonicity bounds for $\mathcal A$ gives
\begin{equation*}
    \|B_kh\|_{V^*}\le L\|h\|_V,
    \qquad
    \langle B_kh,h\rangle_{V^*,V}
    \ge\gamma\|h\|_V^2,
    \qquad h\in V.
\end{equation*}
The same bounds hold for the duality adjoint $B_k^*:V\to V^*$. Given $q\in V^*$, set $h=(B_k^*)^{-1}q$. Since
$\|q\|_{V^*}=\|B_k^*h\|_{V^*}\le L\|h\|_V$,
we obtain the lower bound
\begin{equation*}
    \langle q,(B_k^*)^{-1}q\rangle_{V^*,V}
    =
    \langle B_k^*h,h\rangle_{V^*,V}
    \ge
    \gamma\|h\|_V^2
    \ge
    \frac{\gamma}{L^2}\|q\|_{V^*}^2.
\end{equation*}
For the upper bound, coercivity gives
\begin{equation*}
    \gamma\|h\|_V^2
    \le
    \langle B_k^*h,h\rangle_{V^*,V}
    =
    \langle q,h\rangle_{V^*,V}
    \le
    \|q\|_{V^*}\|h\|_V,
\end{equation*}
so $\|h\|_V\le\gamma^{-1}\|q\|_{V^*}$. Consequently,
\begin{equation*}
    \langle q,(B_k^*)^{-1}q\rangle_{V^*,V}
    \le
    \|q\|_{V^*}\|h\|_V
    \le
    \frac{1}{\gamma}\|q\|_{V^*}^2.
\end{equation*}
Taking $q=\mathcal A(u_k)-f$ proves \eqref{eq:wan_equi_norm}.
\end{proof}

Unlike WAN, NeRI therefore updates only the trial network, without explicitly computing or parameterizing the test function or solving an adversarial maximization problem.

\section{NeRI for elliptic variational inequalities}
\label{sec:VI}

Let $\mathcal A:V\to V^*$ and $f\in V^*$ be defined as in Section~\ref{sec:problem_setting}, and suppose Assumption~\ref{assup:Lip+Coer} holds. For a nonempty, closed, and convex set $K\subset V$, we consider the following elliptic variational inequality.

\begin{definition}\label{def:VI}
Find $u^*\in K$ such that
\begin{equation}\label{eq:VI}
    \langle\mathcal A(u^*)-f,v-u^*\rangle_{V^*,V}
    \ge0,
    \qquad \forall v\in K.
\end{equation}
\end{definition}


We now describe two approaches for incorporating the constraint into NeRI. The first replaces the variational inequality by a penalized operator equation, while the second uses a constraint-preserving neural parameterization.

\subsection{Penalty method}
\label{subsec:VI_penalty}

For concreteness, consider the obstacle set
\begin{equation}\label{eq:obstacle_set}
    K
    :=
    \left\{
        v\in V:\ v\ge\psi
        \text{ a.e. in }\Omega
    \right\},
\end{equation}
where $\psi\in L^2(\Omega)$ is a prescribed obstacle. To penalize violations of the constraint, define
\begin{equation}\label{eq:VI_penalty_functional}
    \mathcal P(u)
    :=
    \frac{1}{2}
    \int_\Omega
    \bigl((\psi-u)^+\bigr)^2\,\mathrm dx,
\end{equation}
where $s^+:=\max\{s,0\}$. Its Fr\'echet derivative satisfies
\begin{equation*}
    \langle D\mathcal P(u),h\rangle_{V^*,V}
    =
    -\int_\Omega(\psi-u)^+h\,\mathrm dx,
    \qquad h\in V.
\end{equation*}
Here and below, $L^2(\Omega)$ terms are interpreted as elements of $V^*$ through the $L^2$ pairing.

For $\varepsilon>0$, define the penalized operator $\mathcal A_\varepsilon(u):= \mathcal A(u)-\frac{1}{\varepsilon}(\psi-u)^+$, and consider the penalized operator equation
\begin{equation}\label{eq:penalized_VI_equation}
    \mathcal A_\varepsilon(u_\varepsilon)=f.
\end{equation}
The map $u\mapsto-(\psi-u)^+$ is monotone and Lipschitz continuous in $L^2(\Omega)$. By the continuous embedding, $\mathcal A_\varepsilon:V\to V^*$ therefore remains strongly monotone and Lipschitz continuous. Under suitable assumptions, standard penalty theory yields $u_\varepsilon\to u^*$ as $\varepsilon\to0$ in an appropriate topology \cite{tremolieres2011numerical,kinderlehrer2000introduction}.

NeRI can now be applied directly to \eqref{eq:penalized_VI_equation}. At iteration $k$, define
\begin{equation}\label{eq:NeRI_VI_loss}
    \mathcal J_{k,\varepsilon}^{\mathrm{VI}}(\theta)
    :=
    \left\langle
        \operatorname{sg}\left(
            \mathcal A(u_{\theta_k})
            -\frac{1}{\varepsilon}(\psi-u_{\theta_k})^+
            -f
        \right),
        u_\theta
    \right\rangle_{V^*,V}.
\end{equation}
The corresponding parameter update is
\begin{equation*}
    \theta_{k+1}
    =
    \theta_k-\eta_k
    \nabla_\theta\mathcal J_{k,\varepsilon}^{\mathrm{VI}}(\theta)
    \big|_{\theta=\theta_k}.
\end{equation*}
Alternatively, the penalty can be included explicitly in the objective:
\begin{equation}\label{eq:loss_VI}
    \widetilde{\mathcal J}_{k,\varepsilon}^{\mathrm{VI}}(\theta)
    :=
    \left\langle
        \operatorname{sg}\bigl(\mathcal A(u_{\theta_k})-f\bigr),
        u_\theta
    \right\rangle_{V^*,V}
    +
    \frac{1}{2\varepsilon}
    \int_\Omega
    \bigl((\psi-u_\theta)^+\bigr)^2\,\mathrm dx.
\end{equation}
Although \eqref{eq:NeRI_VI_loss} and \eqref{eq:loss_VI} differ as functions of $\theta$, their parameter gradients coincide at $\theta_k$. Thus, the penalty contribution can either be included in the detached residual or differentiated explicitly. Decreasing $\varepsilon$ strengthens the penalty but increases the Lipschitz bound for $\mathcal A_\varepsilon$, potentially reducing the admissible step size and worsening the conditioning of the iteration.

\subsection{Hard-constrained parameterization}
\label{subsec:VI_hard_constraint}

We next enforce the obstacle constraint in \eqref{eq:obstacle_set} by construction. If essential boundary conditions are also included in $K$, they must be imposed simultaneously through a suitable boundary-conforming parameterization.

A simple choice is the squared-network parameterization \cite{hintermuller2026constrained}
\begin{equation}\label{eq:hard_VI}
    \widehat u_\theta
    :=
    \psi+\mathcal N_\theta^2.
\end{equation}
Provided that $\psi,\mathcal N_\theta^2\in V$, this ensures $\widehat u_\theta\in K$ for every $\theta$.

For this ansatz, the penalty term in \eqref{eq:loss_VI} vanishes identically, leaving the hard-constrained NeRI objective
\begin{equation}\label{eq:hard_VI_loss}
    \widehat{\mathcal J}_k^{\mathrm{VI}}(\theta)
    :=
    \left\langle
        \operatorname{sg}\bigl(
            \mathcal A(\widehat u_{\theta_k})-f
        \bigr),
        \widehat u_\theta
    \right\rangle_{V^*,V}.
\end{equation}
We then perform the update
\begin{equation*}
    \theta_{k+1}
    =
    \theta_k-\eta_k
    \nabla_\theta\widehat{\mathcal J}_k^{\mathrm{VI}}(\theta)
    \big|_{\theta=\theta_k}.
\end{equation*}
In our numerical experience, the hard-constrained formulation often outperforms the penalty approach, while enforcing feasibility exactly and avoiding the need to tune an obstacle
penalty parameter. Its convergence analysis, however, is more involved, since the nonlinear parameterization changes the tangent geometry and may introduce degeneracies on the contact
set. In particular, stationarity of the parameter iteration does not, in general, imply satisfaction of the original VI. A rigorous analysis of these stationary points and the
convergence of the hard-constrained NeRI iteration is left for future work.

\section{Lazy training and error analysis}
\label{sec:lazy}

In this section, we analyze lazy training for NeRI under gradient flow (depending on an artificial time $t$). We then combine the frozen-kernel approximation with reduced-space estimates to obtain an error decomposition relative to the exact solution.

Let \(\theta_0\) denote the initial network parameters. Writing $G_\theta:=D_\theta u_\theta$ and $K_\theta=G_\theta G_\theta^*$, the continuous-time NeRI dynamics are
\begin{equation}\label{eq:paramflow}
    \dot\theta_t=-G_{\theta_t}^*r_t,
    \qquad
    r_t:=\mathcal A(u_{\theta_t})-f,\qquad
\theta(0)=\theta_0,
\end{equation}
with $\dot\theta_t$ denoting the derivative with respect to $t$.
By the chain rule, $u_t:=u_{\theta_t}$ satisfies the exact function-space equation
\begin{equation}\label{eq:functionflow}
    \dot u_t=-K_{\theta_t}\bigl(\mathcal A(u_t)-f\bigr),\qquad
u(0)=u_{\theta_0}.
\end{equation}

\subsection{An architecture-independent lazy-training theorem}

We first establish a general estimate that isolates the mechanism behind kernel freezing. For a twice Fr\'echet differentiable parameterization $u_\theta:\mathbb R^{p_m}\to V$, define
\begin{equation*}
    H_\theta:=D_\theta^2u_\theta:
    \mathbb R^{p_m}\times\mathbb R^{p_m}\to V,
    \qquad
    \|H_\theta\|_{\mathrm{bil},V}
    :=\sup_{|a|=|b|=1}\|H_\theta[a,b]\|_V.
\end{equation*}
For $R>0$, set
\begin{equation*}
    B_R(\theta_0):=\{\theta:|\theta-\theta_0|\le R\},
    \qquad
    h_{m,R}:=\sup_{\theta\in B_R(\theta_0)}
    \|H_\theta\|_{\mathrm{bil},V}.
\end{equation*}
The following lemma relates Hessian bounds to Jacobian and kernel stability.

\begin{lemma}\label{lem:stability}
Let $\theta,\vartheta\in\mathbb R^{p_m}$. Suppose that, for some $h\ge0,\, \sup_{s\in[0,1]}
\left\|H_{\vartheta+s(\theta-\vartheta)}\right\|_{\mathrm{bil},V}
\le h$. Then
\begin{equation*}
    \|G_\theta-G_\vartheta\|_{\mathcal L(\mathbb R^{p_m},V)}
    \le h|\theta-\vartheta|,\qquad
    \|K_\theta-K_\vartheta\|_{\mathcal L(V^*,V)}
    \le
    \bigl(\|G_\theta\|+\|G_\vartheta\|\bigr)
    h|\theta-\vartheta|.
\end{equation*}
\end{lemma}

\begin{proof}
Fix $q\in\mathbb R^{p_m}$ and put $d:=\theta-\vartheta$. The map $F_q(s):=G_{\vartheta+sd}q$ satisfies
\begin{equation*}
    F_q'(s)=H_{\vartheta+sd}[d,q],
    \qquad
    \|F_q'(s)\|_V\le h|d||q|,
    \qquad 0\le s\le1.
\end{equation*}
The mean-value inequality therefore gives
$\|(G_\theta-G_\vartheta)q\|_V\le h|\theta-\vartheta||q|$.
Taking the supremum over $|q|=1$ proves the first estimate.

For the tangent operator, write
\begin{equation*}
    K_\theta-K_\vartheta
    =(G_\theta-G_\vartheta)G_\theta^*
    +G_\vartheta(G_\theta^*-G_\vartheta^*).
\end{equation*}
The second estimate follows from the first and the identity $\|T^*\|=\|T\|$.
\end{proof}

We next bound the residual and parameter displacement on a finite time interval.

\begin{lemma}\label{lem:path}
Suppose Assumption~\ref{assup:Lip+Coer} holds and $\|G_{\theta_t}\|\le M_G$ for $0\le t\le T$. Then
\begin{align}
    \|r_t\|_{V^*}
    &\le\|r_0\|_{V^*}\exp(LM_G^2t),
    \label{eq:r_t}
    \\
    |\theta_t-\theta_0|
    &\le M_G\|r_0\|_{V^*}
    \int_0^t\exp(LM_G^2s)\,ds.
\end{align}
\end{lemma}

\begin{proof}
Equation~\eqref{eq:functionflow} gives
\begin{equation*}
    \|\dot u_t\|_V
    \le\|K_{\theta_t}\|\|r_t\|_{V^*}
    \le M_G^2\|r_t\|_{V^*}.
\end{equation*}
Integrating and using Lipschitz continuity of $\mathcal A$, we obtain
\begin{equation*}
    \|r_t\|_{V^*}
    \le\|r_0\|_{V^*}+L\|u_t-u_0\|_V
    \le\|r_0\|_{V^*}
    +LM_G^2\int_0^t\|r_s\|_{V^*}\,ds.
\end{equation*}
Since $t\mapsto\|r_t\|_{V^*}$ is continuous, the integral Gronwall inequality yields \eqref{eq:r_t}. The parameter flow then satisfies
\begin{equation*}
    |\dot\theta_t|
    \le\|G_{\theta_t}\|\|r_t\|_{V^*}
    \le M_G\|r_0\|_{V^*}e^{LM_G^2t}.
\end{equation*}
Integrating from $0$ to $t$ proves the displacement bound.
\end{proof}

These estimates allow us to establish lazy training for NeRI. The argument combines Hessian-based stability inside a fixed parameter ball with a bound on the displacement of the trajectory.

\begin{theorem}\label{thm:abstractlazy}
Fix $T>0$ and suppose Assumption~\ref{assup:Lip+Coer} holds. Assume that for every $\delta\in(0,1)$ there exist a deterministic sequence $\varepsilon_m\to0$, events $E_m$ with $\mathbb P(E_m)\ge1-\delta$ for all sufficiently large $m\in\mathbb{N}$, and constants $M_0,C_0,R_0>0$ independent of $m$, such that on $E_m$
\begin{equation*}
    \|G_{\theta_0}\|\le M_0,
    \qquad
    \|r_0\|_{V^*}\le C_0,
    \qquad
    h_{m,R_0}\le\varepsilon_m.
\end{equation*}
If
\begin{equation*}
    R_0>R_T:=2M_0C_0\int_0^T\exp(4LM_0^2s)\,ds,
\end{equation*}
then, on $E_m$ for all sufficiently large $m$, the trajectory exists on $[0,T]$, remains in $B_{R_0}(\theta_0)$, and satisfies
\begin{equation}\label{eq:kernel-drift-general}
    \sup_{0\le t\le T}
    \|K_{\theta_t}-K_{\theta_0}\|_{\mathcal L(V^*,V)}
    \le3M_0R_T\varepsilon_m.
\end{equation}
In particular, the tangent-operator drift converges to zero in probability.
\end{theorem}

\begin{proof}
Fix an event $E_m$ on which the assumed bounds hold and choose $m$ large enough that $\varepsilon_mR_0\le M_0$. By Lemma~\ref{lem:stability}, the vector field in \eqref{eq:paramflow} is locally Lipschitz inside $B_{R_0}(\theta_0)$, so a unique local solution exists.

Let $\tau$ be the first time the trajectory reaches the boundary of this ball, with $\tau=+\infty$ if no such time exists. Before this exit, Lemma~\ref{lem:stability} gives
\begin{equation*}
    \|G_{\theta_t}-G_{\theta_0}\|
    \le\varepsilon_m|\theta_t-\theta_0|
    \le\varepsilon_mR_0\le M_0,
\end{equation*}
and hence $\|G_{\theta_t}\|\le2M_0$. Applying Lemma~\ref{lem:path} with $M_G=2M_0$ yields, for $t\le T$ before the exit,
\begin{equation*}
    |\theta_t-\theta_0|
    \le2M_0C_0\int_0^t e^{4LM_0^2s}\,ds
    \le R_T<R_0.
\end{equation*}
At an exit time $\tau\le T$, this would contradict $|\theta_\tau-\theta_0|=R_0$. The trajectory thus stays in a compact subset of the ball and extends to $T$. In particular,
\begin{equation*}
    \sup_{0\le t\le T}|\theta_t-\theta_0|\le R_T,
    \qquad
    \sup_{0\le t\le T}\|G_{\theta_t}\|\le2M_0.
\end{equation*}
A second application of Lemma~\ref{lem:stability} gives
\begin{equation*}
    \|K_{\theta_t}-K_{\theta_0}\|
    \le(\|G_{\theta_t}\|+\|G_{\theta_0}\|)
    \varepsilon_m|\theta_t-\theta_0|
    \le3M_0R_T\varepsilon_m.
\end{equation*}
This proves \eqref{eq:kernel-drift-general}. Since $\delta>0$ is arbitrary and $\varepsilon_m\to0$, the drift converges to zero in probability.
\end{proof}

We next compare NeRI with the frozen-kernel flow. Set $K_m^0:=K_{\theta_0}$ and let $v_t$ solve
\begin{equation}\label{eq:frozenflow}
    \dot v_t=-K_m^0\bigl(\mathcal A(v_t)-f\bigr),
    \qquad v_0=u_{\theta_0}.
\end{equation}
To relate this flow to the exact solution, define
\begin{equation*}
    X_m:=\operatorname{Ran}(K_m^0),
    \qquad
    S_m:=u_{\theta_0}+X_m,
    \qquad
    B_m:=K_m^0\mathcal R^{-1}:V\to V.
\end{equation*}
The operator $B_m$ is self-adjoint and positive definite on its finite-dimensional range $X_m$. Hence there exist $0<\gamma_{K,m}\le M_{K,m}$ such that
\begin{equation}\label{eq:kernel-spectrum}
    \gamma_{K,m}\|w\|_V^2
    \le\langle B_mw,w\rangle_V
    \le M_{K,m}\|w\|_V^2,
    \qquad w\in X_m.
\end{equation}
Let $u_{S_m}^*\in S_m$ denote the unique reduced solution satisfying
\begin{equation}\label{eq:reduced-frozen}
    \langle\mathcal A(u_{S_m}^*)-f,w\rangle_{V^*,V}=0,
    \qquad w\in X_m.
\end{equation}
The Moore--Penrose pseudoinverse $B_m^\dagger$ defines the norm
\begin{equation*}
    \|w\|_{B_m^\dagger}^2
    :=\langle B_m^\dagger w,w\rangle_V,
    \qquad w\in X_m.
\end{equation*}

\begin{theorem}\label{thm:frozen}
Under the assumptions of Theorem~\ref{thm:abstractlazy}, there is a constant $C_T$, independent of $m$, such that on the same high-probability event, for all sufficiently large $m$,
\begin{equation}\label{eq:frozen-approx}
    \sup_{0\le t\le T}\|u_{\theta_t}-v_t\|_V
    \le C_T\varepsilon_m.
\end{equation}
Moreover, the frozen flow satisfies, for $t\ge0$,
\begin{equation}\label{eq:frozen-reduced-error}
    \|v_t-u_{S_m}^*\|_{B_m^\dagger}
    \le e^{-\gamma\gamma_{K,m}t}
    \|v_0-u_{S_m}^*\|_{B_m^\dagger}.
\end{equation}
If $u^*$ solves $\mathcal A(u^*)=f$, then
\begin{equation}\label{eq:frozen-exact-error}
    \|v_t-u^*\|_V
    \le\sqrt{M_{K,m}}\,e^{-\gamma\gamma_{K,m}t}
    \|v_0-u_{S_m}^*\|_{B_m^\dagger}
    +\frac{L}{\gamma}\inf_{w\in S_m}\|w-u^*\|_V.
\end{equation}
Consequently, for $0\le t\le T$,
\begin{equation}\label{eq:total-error}
    \|u_{\theta_t}-u^*\|_V
    \le C_T\varepsilon_m
    +\sqrt{M_{K,m}}\,e^{-\gamma\gamma_{K,m}t}
    \|u_{\theta_0}-u_{S_m}^*\|_{B_m^\dagger}
    +\frac{L}{\gamma}\inf_{w\in S_m}\|w-u^*\|_V.
\end{equation}
\end{theorem}

\begin{proof}
The frozen equation is globally well posed because $\mathcal A$ is globally Lipschitz. Set $\xi_t:=u_{\theta_t}-v_t$. Subtracting the evolution equations gives
\begin{equation*}
    \dot\xi_t
    =-(K_{\theta_t}-K_m^0)r_t
    -K_m^0\bigl(\mathcal A(u_{\theta_t})-\mathcal A(v_t)\bigr),
    \qquad\xi_0=0.
\end{equation*}
Hence,
\begin{equation*}
    \|\dot\xi_t\|_V
    \le\|K_{\theta_t}-K_m^0\|\|r_t\|_{V^*}
    +L\|K_m^0\|\|\xi_t\|_V.
\end{equation*}
On the event of Theorem~\ref{thm:abstractlazy},
\begin{equation*}
    \sup_{t\le T}\|K_{\theta_t}-K_m^0\|
    \le3M_0R_T\varepsilon_m,
    \qquad
    \|K_m^0\|\le M_0^2,
\end{equation*}
and Lemma~\ref{lem:path} gives
$\sup_{t\le T}\|r_t\|_{V^*}\le C_0e^{4LM_0^2T}$.
Thus, for a constant $A_T$ independent of $m$,
\begin{equation*}
    \|\xi_t\|_V
    \le A_T\varepsilon_m t
    +LM_0^2\int_0^t\|\xi_s\|_V\,ds.
\end{equation*}
The integral Gronwall inequality \cite{Amann1990} yields $\|\xi_t\|_V\le A_TT e^{LM_0^2T}\varepsilon_m $ for $0\le t\le T,$ proving \eqref{eq:frozen-approx}.

We next estimate the frozen flow itself. Since $\dot v_t\in X_m$, we have $v_t\in S_m$ for all $t\ge0$. Set $e_t:=v_t-u_{S_m}^*\in X_m$. By \eqref{eq:reduced-frozen},
$K_m^0\bigl(\mathcal A(u_{S_m}^*)-f\bigr)=0$, and hence
\begin{equation*}
    \dot e_t=-K_m^0\bigl(\mathcal A(v_t)-\mathcal A(u_{S_m}^*)\bigr).
\end{equation*}
Using $K_m^0=B_m\mathcal R$ and $e_t\in X_m$, we obtain
\begin{align*}
    \frac12\frac{d}{dt}\|e_t\|_{B_m^\dagger}^2
    &=-\langle\mathcal A(v_t)-\mathcal A(u_{S_m}^*),e_t\rangle_{V^*,V}\le-\gamma\|e_t\|_V^2
    \le-\gamma\gamma_{K,m}\|e_t\|_{B_m^\dagger}^2.
\end{align*}
Gronwall's inequality proves \eqref{eq:frozen-reduced-error}, while \eqref{eq:kernel-spectrum} gives
$\|e_t\|_V\le\sqrt{M_{K,m}}\,\|e_t\|_{B_m^\dagger}$.

Finally, applying the reduced-space estimate \eqref{eq:finite_rank_cea} with $S=S_m$ yields
\begin{equation}\label{eq:reduced-approx}
    \|u_{S_m}^*-u^*\|_V
    \le\frac{L}{\gamma}\inf_{w\in S_m}\|w-u^*\|_V.
\end{equation}
Combining this bound with \eqref{eq:frozen-reduced-error} proves \eqref{eq:frozen-exact-error}. The total error estimate \eqref{eq:total-error} follows from \eqref{eq:frozen-approx} and the triangle inequality.
\end{proof}

Estimate~\eqref{eq:total-error} separates the total error into finite-width lazy-training, frozen-flow optimization, and reduced-space approximation errors.

\subsection{Two-layer networks in \texorpdfstring{$H^1(\Omega)$}{H1(Omega)}}

We now verify the abstract hypotheses for a two-layer network in $V=H^1(\Omega)$, where $\Omega\subset\mathbb R^d$ is bounded and $d$ is fixed. Consider the two-layer network
\begin{equation}\label{eq:two-layer}
    u_\theta^{(m)}(x)
    =\frac{1}{\sqrt m}\sum_{i=1}^m
    a_i\,\sigma(w_i\cdot x+b_i).
\end{equation}

We use the stochastic order notation of van der Vaart \cite{vanDerVaart1998}. For random elements $Z_m$ of a normed space $(\mathcal X,\|\cdot\|_{\mathcal X})$, we write $Z_m=O_{\mathbb P}(1)$ if they are bounded in probability: for every $\varepsilon>0$, there exist $M_\varepsilon<\infty$ and $m_\varepsilon$ such that $\mathbb P\bigl(\|Z_m\|_{\mathcal X}>M_\varepsilon\bigr)
\le\varepsilon$ for all $m\ge m_\varepsilon$. For a deterministic sequence $a_m>0$, the notation $Z_m=O_{\mathbb P}(a_m)$ means $Z_m/a_m=O_{\mathbb P}(1)$, while $Z_m=o_{\mathbb P}(1)$ means $\|Z_m\|_{\mathcal X}\to0$ in probability.

\begin{assumption}\label{ass:two-layer}
The activation satisfies $\sigma\in C_b^3(\mathbb R)$. At initialization, the blocks $(a_i,w_i,b_i)$ are independent and identically distributed, all coordinates are sub-Gaussian with width-independent sub-Gaussian norms, and $a_i$ is centered and independent of $(w_i,b_i)$.
\end{assumption}

The centering and independence of the output weights ensure that the initial network is $O_{\mathbb P}(1)$ in $H^1(\Omega)$, as shown below.

\begin{lemma}\label{lem:init-two}
Under Assumptions~\ref{ass:two-layer} and~\ref{assup:Lip+Coer},
\begin{equation*}
    \|u_{\theta_0}^{(m)}\|_{H^1}=O_{\mathbb P}(1),
    \qquad
    \|G_{\theta_0}^{(m)}\|_{\mathcal L(\mathbb R^{p_m},H^1)}
    =O_{\mathbb P}(1),
    \qquad
    \|\mathcal A(u_{\theta_0}^{(m)})-f\|_{H^1(\Omega)^*}
    =O_{\mathbb P}(1).
\end{equation*}
\end{lemma}

\begin{proof}
Write $z_i(x)=w_i\cdot x+b_i$. By centering and independence, all cross terms vanish in the second moment of the width-normalized sum. Hence,
\begin{align*}
    \mathbb E\|u_{\theta_0}^{(m)}\|_{L^2(\Omega)}^2
    &=\frac1m\sum_{i=1}^m
    \mathbb E\int_\Omega a_i^2\sigma(z_i(x))^2\,dx
    \le C|\Omega|,
    \\
    \mathbb E\|\nabla u_{\theta_0}^{(m)}\|_{L^2(\Omega)}^2
    &=\frac1m\sum_{i=1}^m
    \mathbb E\int_\Omega a_i^2\sigma'(z_i(x))^2|w_i|^2\,dx
    \le C.
\end{align*}
The relevant initialization moments are uniformly bounded, so $C$ is independent of $m$. Markov's inequality gives the first estimate.

For the Jacobian, the parameter derivatives are
\begin{equation*}
    \partial_{a_i}u=m^{-1/2}\sigma(z_i),
    \qquad
    \partial_{w_{i,j}}u=m^{-1/2}a_i\sigma'(z_i)x_j,
    \qquad
    \partial_{b_i}u=m^{-1/2}a_i\sigma'(z_i).
\end{equation*}
Their spatial gradients are
\begin{align*}
    \nabla\partial_{a_i}u
    &=m^{-1/2}\sigma'(z_i)w_i,
    \qquad
    \nabla\partial_{b_i}u=m^{-1/2}a_i\sigma''(z_i)w_i,
    \\
    \nabla\partial_{w_{i,j}}u
    &=m^{-1/2}a_i
    \bigl[\sigma''(z_i)x_jw_i+\sigma'(z_i)e_j\bigr].
\end{align*}
Since $\Omega$ and $\sigma,\sigma',\sigma''$ are bounded, summing the squared $H^1$ norms gives
\begin{equation*}
    \|G_{\theta_0}^{(m)}\|_{\mathrm{HS}(\mathbb R^{p_m},H^1)}^2
    \le\frac{C}{m}\sum_{i=1}^m
    (1+a_i^2)(1+|w_i|^2).
\end{equation*}
Uniform moment bounds and Markov's inequality make the right-hand side $O_{\mathbb P}(1)$. Since the operator norm is bounded by the Hilbert--Schmidt norm, the Jacobian estimate follows. Finally, the residual bound follows from the following inequality:
\begin{equation*}
    \|\mathcal A(u_{\theta_0}^{(m)})-f\|_{H^1(\Omega)^*}
    \le L\|u_{\theta_0}^{(m)}\|_{H^1}
    +\|\mathcal A(0)-f\|_{H^1(\Omega)^*}.
\end{equation*}
\end{proof}

We now derive the following Hessian estimate.

\begin{lemma}\label{lem:hessian-two}
For fixed $R>0$, Assumption~\ref{ass:two-layer} implies
\begin{equation}\label{eq:hessian-rate}
    \sup_{\theta\in B_R(\theta_0)}
    \|D_\theta^2u_\theta^{(m)}\|_{\mathrm{bil},H^1}
    =O_{\mathbb P}\!\left(\frac{\log m}{\sqrt m}\right).
\end{equation}
\end{lemma}

\begin{proof}
For one neuron, write
\begin{equation*}
    \phi(a,w,b;x)=m^{-1/2}a\sigma(z),
    \qquad z=w\cdot x+b.
\end{equation*}
For block directions $p=(\alpha,r,\beta)$ and $q=(\widetilde\alpha,\widetilde r,\widetilde\beta)$, set
$\rho_p=r\cdot x+\beta$ and $\rho_q=\widetilde r\cdot x+\widetilde\beta$.
Direct differentiation gives
\begin{equation}\label{eq:two-neuron-hess}
    D^2\phi[p,q]
    =\frac1{\sqrt m}
    \bigl[
        \alpha\sigma'(z)\rho_q
        +\widetilde\alpha\sigma'(z)\rho_p
        +a\sigma''(z)\rho_p\rho_q
    \bigr].
\end{equation}
Boundedness of $\Omega$ gives $|\rho_p(x)|\le C_\Omega|p|$ and $|\rho_q(x)|\le C_\Omega|q|$, so
\begin{equation*}
    \|D^2\phi[p,q]\|_{L^2(\Omega)}
    \le\frac{C}{\sqrt m}(1+|a|)|p||q|.
\end{equation*}
Taking one spatial derivative yields
\begin{align*}
    \nabla_xD^2\phi[p,q]
    =\frac1{\sqrt m}\Bigl(&
        \alpha[\sigma''(z)w\rho_q+\sigma'(z)\widetilde r]
        +\widetilde\alpha[\sigma''(z)w\rho_p+\sigma'(z)r]
        \\
        &+a[\sigma'''(z)w\rho_p\rho_q
        +\sigma''(z)r\rho_q
        +\sigma''(z)\widetilde r\rho_p]
    \Bigr).
\end{align*}
Since $\sigma',\sigma'',\sigma'''$ are bounded,
\begin{equation*}
    \|\nabla_xD^2\phi[p,q]\|_{L^2(\Omega)}
    \le\frac{C}{\sqrt m}(1+|a|)(1+|w|)|p||q|.
\end{equation*}

Partition global directions as $p=(p_1,\ldots,p_m)$ and $q=(q_1,\ldots,q_m)$. Since cross-neuron derivatives vanish,
$D_\theta^2u_\theta[p,q]=\sum_{i=1}^mD^2\phi_i[p_i,q_i]$.
The preceding one-neuron bounds and the triangle inequality give
\begin{align*}
    \|D_\theta^2u_\theta[p,q]\|_{H^1}
    &\le\frac{C}{\sqrt m}
    \max_{i\le m}(1+|a_i|)(1+|w_i|)
    \sum_{i=1}^m|p_i||q_i|
    \\
    &\le\frac{C}{\sqrt m}
    \max_{i\le m}(1+|a_i|)(1+|w_i|)|p||q|,
\end{align*}
where the last step is Cauchy--Schwarz over the neuron blocks. For $\theta\in B_R(\theta_0)$,
\begin{equation*}
    |a_i|\le|a_i(0)|+R,
    \qquad
    |w_i|\le|w_i(0)|+R.
\end{equation*}
The sub-Gaussian maximal bounds, with $d$ fixed, yield
\begin{equation*}
    \max_{i\le m}|a_i(0)|=O_{\mathbb P}(\sqrt{\log m}),
    \qquad
    \max_{i\le m}|w_i(0)|=O_{\mathbb P}(\sqrt{\log m}).
\end{equation*}
Thus the product factor is $O_{\mathbb P}(\log m)$ uniformly on $B_R(\theta_0)$. Taking the supremum over unit $p,q$ proves \eqref{eq:hessian-rate}.
\end{proof}

\begin{theorem}\label{thm:two-layer}
Let $\Omega\subset\mathbb R^d$ be bounded and suppose Assumptions~\ref{assup:Lip+Coer} and~\ref{ass:two-layer} hold. For every fixed $T>0$,
\begin{equation*}
    \sup_{0\le t\le T}
    \|K_{\theta_t}^{(m)}-K_{\theta_0}^{(m)}\|_{\mathcal L(H^1(\Omega)^*,H^1(\Omega))}
    =O_{\mathbb P}\!\left(\frac{\log m}{\sqrt m}\right).
\end{equation*}
If $v_t^{(m)}$ solves \eqref{eq:frozenflow}, then
\begin{equation*}
    \sup_{0\le t\le T}
    \|u_{\theta_t}^{(m)}-v_t^{(m)}\|_{H^1}
    =O_{\mathbb P}\!\left(\frac{\log m}{\sqrt m}\right).
\end{equation*}
Consequently, for $0\le t\le T$,
\begin{align}
    \|u_{\theta_t}^{(m)}-u^*\|_{H^1}
    \le{}&\sqrt{M_{K,m}}\,e^{-\gamma\gamma_{K,m}t}
    \|u_{\theta_0}^{(m)}-u_{S_m}^*\|_{B_m^\dagger}
    \nonumber\\
    &+\frac{L}{\gamma}\inf_{w\in S_m}\|w-u^*\|_{H^1}
    +O_{\mathbb P}\!\left(\frac{\log m}{\sqrt m}\right),
    \label{eq:two-layer-total-error}
\end{align}
where the stochastic remainder is uniform on $[0,T]$.
\end{theorem}

\begin{proof}
Fix $T>0$ and $\delta>0$. Lemma~\ref{lem:init-two} gives constants $M_0,C_0>0$, independent of $m$, such that
\begin{equation*}
    \mathbb P\!\left(
        \|G_{\theta_0}^{(m)}\|\le M_0,
        \quad
        \|\mathcal A(u_{\theta_0}^{(m)})-f\|_{H^1(\Omega)^*}\le C_0
    \right)
    \ge1-\delta/2
\end{equation*}
for all sufficiently large $m$. Choose $R_0>R_T$, with $R_T$ defined from these constants in Theorem~\ref{thm:abstractlazy}. Applying Lemma~\ref{lem:hessian-two} on $B_{R_0}(\theta_0)$ gives a constant $C_\delta$ such that
\begin{equation*}
    \mathbb P\!\left(
        \sup_{\theta\in B_{R_0}(\theta_0)}
        \|D_\theta^2u_\theta^{(m)}\|_{\mathrm{bil},H^1}
        \le C_\delta\frac{\log m}{\sqrt m}
    \right)
    \ge1-\delta/2.
\end{equation*}
The intersection has probability at least $1-\delta$, and Theorem~\ref{thm:abstractlazy} applies with $\varepsilon_m=C_\delta(\log m)/\sqrt m$. The kernel-drift bound follows from \eqref{eq:kernel-drift-general}, and Theorem~\ref{thm:frozen} gives the frozen-flow approximation. Since $\delta>0$ is arbitrary, both estimates have the stated stochastic order. Finally, \eqref{eq:two-layer-total-error} follows from \eqref{eq:total-error}.
\end{proof}

Estimate~\eqref{eq:two-layer-total-error} separates the total error into finite-width lazy-training, frozen-flow optimization, and reduced-space approximation errors. The frozen-flow approximation holds for fixed $T$, with constants that may depend on $T$. These estimates do not by themselves establish convergence to $u^*$ in a joint large-width and long-time limit. Such convergence requires further control of the reduced-space approximation error and compatibility between the width-dependent optimization time and the validity interval of the frozen-flow approximation.

\section{Numerical results}
\label{sec:numerical}

In this section, we evaluate the performance of NeRI on several representative problems and compare it with PINN and WAN. Unless otherwise specified, all experiments use a residual neural network architecture consisting of $N_{\mathrm{blk}}$ residual blocks, each containing two fully connected layers and a residual connection. We set $N_{\mathrm{blk}}=2$ for two-dimensional problems and $N_{\mathrm{blk}}=3$ for high-dimensional problems. Each hidden layer has width $100$, and the hyperbolic tangent function is used as the activation throughout.

All networks are trained using Adam, and we use cosine annealing to decrease the learning rate from $1.0\times 10^{-3}$ to $1.0\times 10^{-8}$ over the full training run. At each training epoch, the interior and boundary points are sampled independently and uniformly in $\Omega$ and on $\partial\Omega$, respectively. We use $N_{\Omega}=10{,}000$ interior points in all experiments, with 256 points per boundary segment (including the slit, when present) in two dimensions and 2,000 points per boundary face in higher dimensions. The experiments were implemented in PyTorch and executed on a MacBook Pro equipped with an Apple M5 chip. The implementation code of the numerical experiments is publicly available at \url{https://github.com/JianfengNing/NeRI-Code}.

\subsection{A nonsymmetric elliptic problem with a corner singularity}
\label{nume:corner}

Let $\Omega=(-1,1)^2\setminus([0,1)\times\{0\})$ be the slit domain, consider the nonsymmetric elliptic problem
\begin{equation*}
    \begin{cases}
        -\Delta u+u_{x_1}=f & \text{in } \Omega,\\
        u=g & \text{on } \partial\Omega.
    \end{cases}
\end{equation*}
The tip of the slit at the origin gives rise to a corner singularity. We prescribe the exact solution $u^*(r,\theta)=r^{1/2}\sin\left(\frac{\theta}{2}\right)$, where $(r,\theta)$ are polar coordinates centered at the origin, and choose $f$ and $g$ accordingly. The solution belongs to $H^1(\Omega)$ but not to $H^2(\Omega)$. The crack-tip singularity can pose approximation and optimization challenges for standard strong-form PINNs, whereas NeRI and WAN use weak formulations involving only first-order spatial derivatives.

We compare NeRI with PINN and WAN. For all three methods, the Dirichlet boundary condition is imposed through a quadratic penalty with parameter $\beta/2=50$, and the networks are trained for $15{,}000$ epochs. At each epoch, WAN performs one minimization update and two maximization updates.

\begin{figure}[htbp]
    \centering

    \begin{subfigure}[b]{0.8\textwidth}
        \centering
        \includegraphics[width=\linewidth]{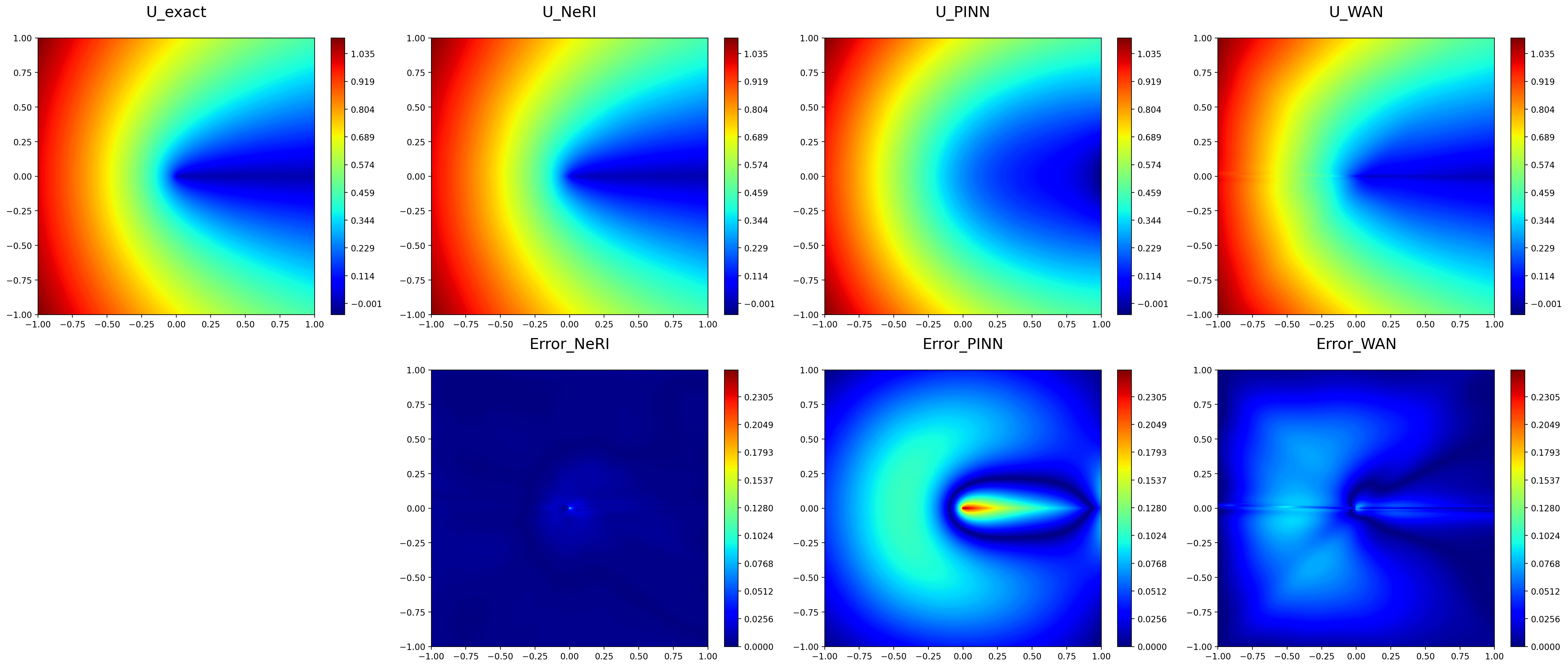}
        \caption{The first row presents the exact solution, numerical solutions obtained by NeRI, PINN, and WAN. The second row presents the corresponding absolute errors.}
        \label{fig:corner}
    \end{subfigure}
    \begin{subfigure}[b]{0.5\textwidth}
        \centering
        \includegraphics[width=\linewidth]{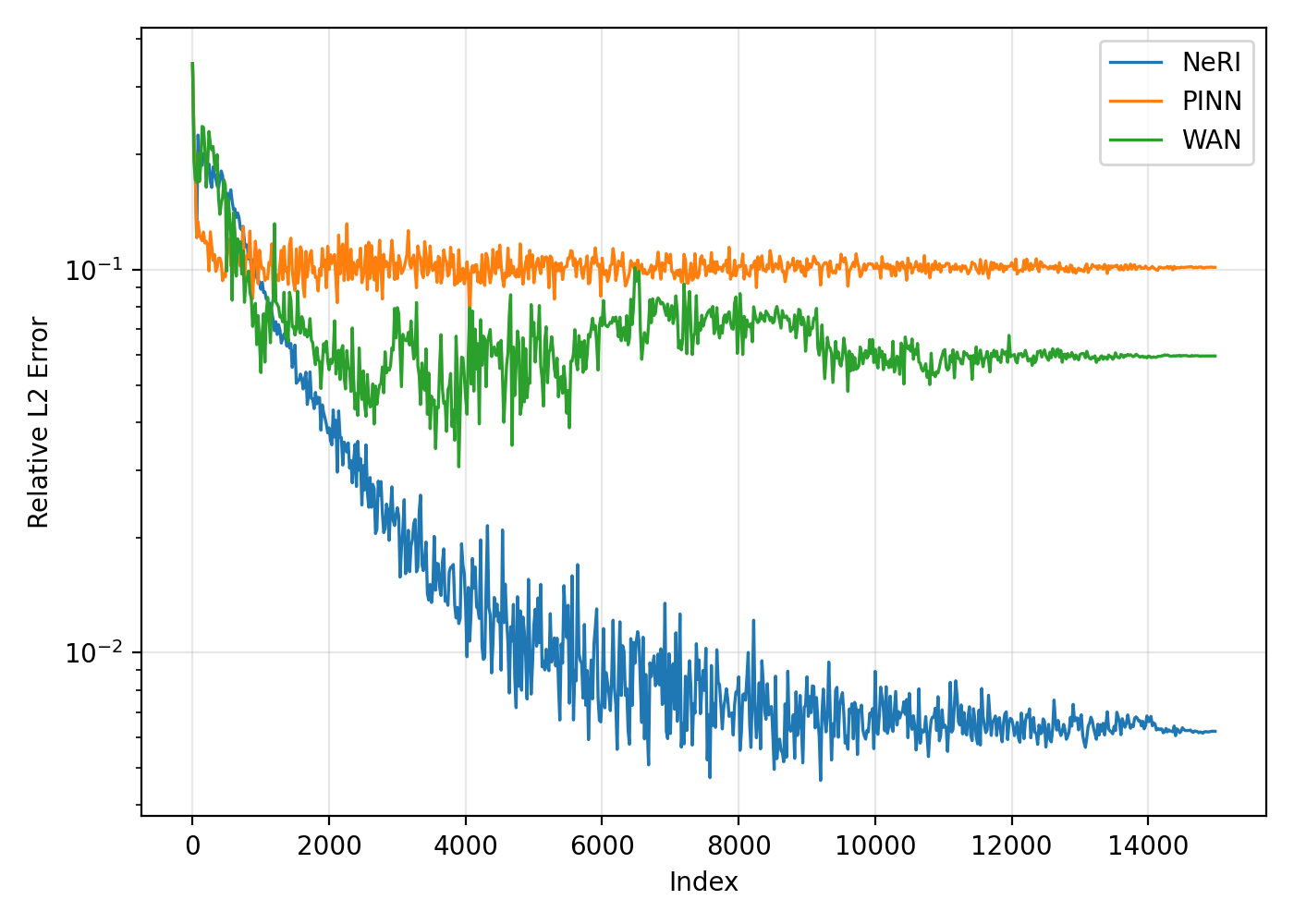}
        \caption{Relative $L^2$ error histories during training.}
        \label{fig:error_history_corner}
    \end{subfigure}

    \caption{Numerical results for Example \ref{nume:corner}.}
    \label{fig:corner_results}
\end{figure}

The numerical approximations and corresponding pointwise absolute errors are shown in Fig.~\ref{fig:corner}. NeRI yields the smallest error among the three methods, while PINN exhibits noticeably larger errors near the singular point and WAN gives an intermediate accuracy. The relative $L^2$ errors obtained by NeRI, PINN, and WAN are
$6.24\times10^{-3}$, $1.01\times10^{-1}$, and $5.95\times10^{-2}$, respectively. Fig.~\ref{fig:error_history_corner} shows the corresponding error histories during training. The error of NeRI decreases steadily and remains stable throughout the optimization. PINN decreases rapidly at the initial stage but subsequently stagnates at a comparatively high error level, consistent with the difficulty of strong-form training for solutions with limited regularity. WAN also shows a stable initial decrease, but reaches a higher error level than NeRI, illustrating the additional optimization difficulty associated with its min--max formulation.

To more directly connect the numerical experiments with the theoretical setting of Section~\ref{sec:lazy}, we additionally test NeRI using plain gradient descent and two-layer networks with different widths. The learning rate is fixed at $5\times10^{-3}$, and all networks are trained for $100{,}000$ epochs. As is commonly observed in neural PDE solvers, adaptive optimizers such as Adam often exhibit more favorable empirical optimization behavior than vanilla gradient descent. Consistently, Table~\ref{tab:two_layer} shows that plain gradient descent leads to a larger error than Adam for NeRI in this example, while still achieving reasonable accuracy. The resulting errors also remain below the PINN and WAN errors reported above, although the corresponding architectures and optimizers are different.

\begin{table}[htbp]
 \caption{Relative $L^2$ errors of NeRI for Example~\ref{nume:corner}
using plain gradient descent and two-layer networks with different widths.}
 \label{tab:two_layer}
 \centering
 \small
 \setlength{\tabcolsep}{5pt}
 \begin{tabular}{@{}ccccc@{}}
 \toprule
 Width & 25 & 50 & 100 & 200 \\
 \midrule
Relative $L^2$ Error & $4.84\times 10^{-2}$ & $4.43\times 10^{-2}$  & $4.09\times 10^{-2}$ &$4.23\times 10^{-2}$\\
 \bottomrule
 \end{tabular}
\end{table}

\subsection{A low-regularity problem with a discontinuous diffusion coefficient}
\label{nume:diffusion}

Let $\Omega=(0,1)^2$ and consider the nonsymmetric elliptic problem
\begin{equation*}
\begin{cases}
\begin{aligned}
    -\nabla\cdot\bigl(a(x)\nabla u\bigr)+u_{x_1}
    &=10 \quad \text{in }\Omega,\\
    u&=0 \quad \text{on }\partial\Omega.
\end{aligned}
\end{cases}
\end{equation*}
The diffusion coefficient is piecewise constant with $a(x)=3$ for $\frac{1}{3}\le x_1\le\frac{2}{3}$, and $a(x)=1$ otherwise. The jumps in $a$ across the interfaces reduce the global regularity of the solution and make the problem challenging for standard strong-form neural solvers, whereas the weak formulation naturally accommodates the discontinuous coefficient.

\begin{figure}
    \centering
    \includegraphics[width=0.8\linewidth]{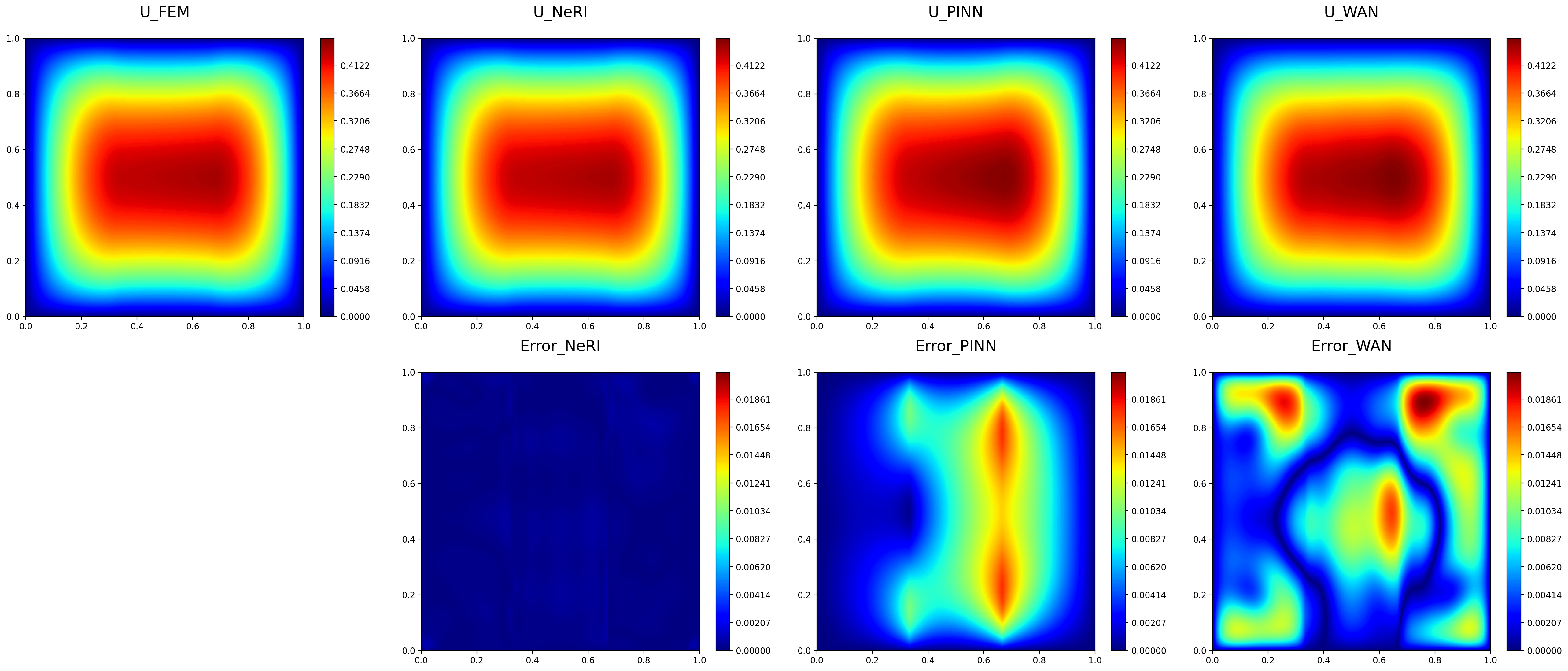}
    \caption{Numerical results for Example~\ref{nume:diffusion}
    . The first row presents the exact solution, numerical approximations obtained by NeRI, PINN, and WAN. The second row presents the corresponding absolute errors.}
    \label{fig:Diffusion}
\end{figure}

We compare NeRI with PINN and WAN. For all three methods, the homogeneous Dirichlet condition is imposed through the hard-constrained parameterization $u_\theta=\mathcal{N}_\theta w$ with $w(x)=x_1(1-x_1)x_2(1-x_2)$. The networks are trained for $20{,}000$ epochs. The numerical approximations and corresponding pointwise absolute errors are shown in Fig.~\ref{fig:Diffusion}. The reference solution is computed using $P_1$ finite elements on a triangulated uniform $300\times300$ grid with centroid quadrature. NeRI gives the smallest error among the three methods and captures the solution behavior across the coefficient interfaces more accurately. PINN exhibits noticeably larger errors, especially near the interfaces. The relative $L^2$ errors for NeRI, PINN, and WAN are
$9.13\times10^{-4}$, $2.70\times10^{-2}$, and $2.89\times10^{-2}$, respectively.

We emphasize that the PINN used here is a standard strong-form, single-network baseline. In an additional experiment, an interface-aware PINN using three subnetworks and penalizing jumps in the solution and the normal diffusive flux achieved a relative \(L^2\) error of \(1.63\times10^{-3}\). Such approaches, however, require explicit treatment and sampling of the interfaces and may become more involved for complex geometries or high-dimensional problems. In contrast, NeRI and WAN incorporate the discontinuous coefficient directly through the weak formulation without requiring explicit interface conditions.

\subsection{High-dimensional $p$-Laplace problem}
\label{nume:p_laplace}

\begin{figure}
    \centering
    \includegraphics[width=0.75\linewidth]{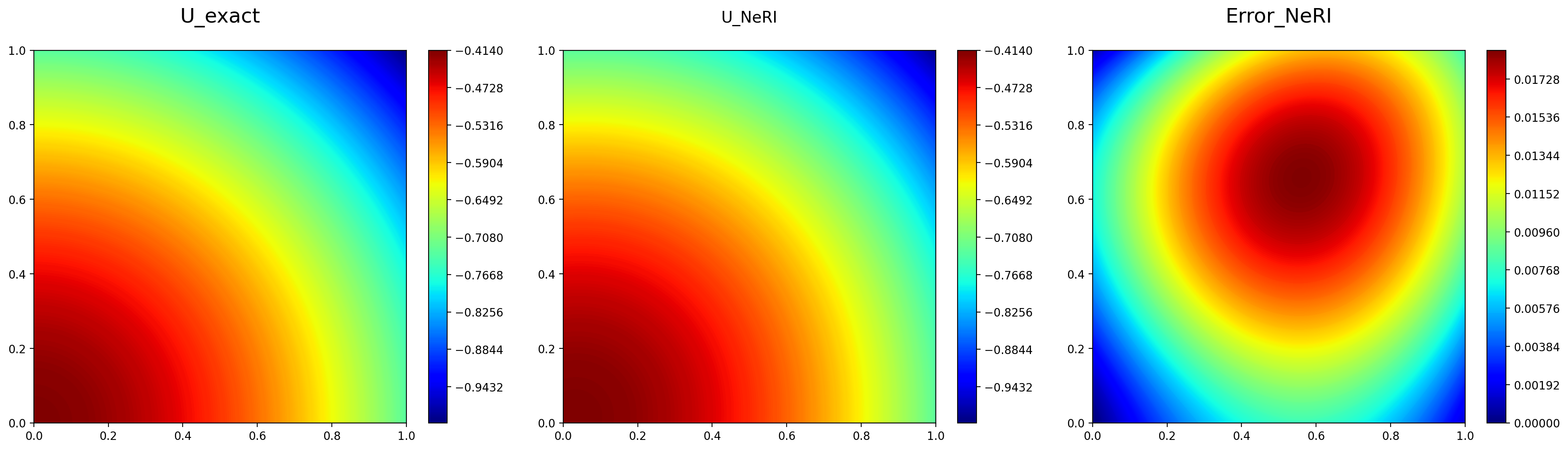}
    \caption{Numerical results for Example \ref{nume:p_laplace}.}
    \label{fig:p_laplace}
\end{figure}

Let $\Omega=(0,1)^{10}$. We consider the following $p$-Laplace problem studied in \cite{hu2025iterative}:
\begin{equation*}
    \begin{cases}
    \begin{aligned}
        -\nabla\cdot\bigl(|\nabla u|^{p-2}\nabla u\bigr)
        +\vec b\cdot\nabla u
        &=f \quad \text{in }\Omega,\\
        u&=g \quad \text{on }\partial\Omega.
    \end{aligned}
    \end{cases}
\end{equation*}
We set $p=1.5$, $f(x)=\frac{9-x_1}{\|x\|_2}$ and $\vec{b}=(1,0,\cdots,0)$. The exact solution is given as $u^*(x)=1-\|x\|_2$. The network is trained for $6{,}000$ epochs, and the Dirichlet boundary condition is imposed through a quadratic penalty with parameter $\beta=200$. The numerical solution on the two-dimensional slice
$(x_1,x_2,0.5,\ldots,0.5)$ is shown in Fig.~\ref{fig:p_laplace}.

The relative $L^2$ error of NeRI over $\Omega$ is $8.04\times10^{-3}$, and the relative $L^2$ error on the same slice is $1.95\times10^{-2}$. For this setting, \cite{hu2025iterative} reports errors of $2.86\times10^{-2}$, $6.71\times10^{-2}$, and $1.03\times10^{-1}$ for IDRM, PINN, and WAN, respectively. Moreover, compared with IDRM, NeRI avoids inner surrogate minimization and the associated regularization and stopping parameters. It updates the network directly from the current weak residual, without differentiating through the nonlinear operator, yielding a simpler training procedure.

We note that the \(p\)-Laplace operator with \(p=1.5\) does not satisfy Assumption~\ref{assup:Lip+Coer}, so this example lies outside the error analysis of Section~\ref{sec:lazy}. Nevertheless, the weak adversarial interpretation remains valid where the implicit adjoint problem is well posed and the linearized operator has a strictly positive quadratic form: the pairing is then strictly positive for nonzero residuals, even without uniform equivalence to the squared dual residual norm. Numerical results indicate that NeRI remains effective for the tested problems, while extending the analysis to this broader setting is left for future work.

\subsection{High-dimensional biharmonic equation}
\label{nume:biharmonic}

Let $\Omega=(0,1)^{10}$ and consider the following high-dimensional biharmonic problem:
\begin{equation*}
    \begin{cases}
    \begin{aligned}
        \Delta^2 u+\vec b\cdot\nabla u&=f
        \quad \text{in }\Omega,\\
        \Delta u=0,\qquad u&=0
        \quad \text{on }\partial\Omega.
    \end{aligned}
    \end{cases}
\end{equation*}
We take $\vec b=(1,1,0,\ldots,0)$ and set
$V=H^2(\Omega)\cap H_0^1(\Omega)$. The corresponding weak formulation is
\begin{equation}
    \langle\mathcal A u-f,v\rangle_{V^*,V}
    =
    \int_\Omega \Delta u\,\Delta v\,\mathrm dx
    +\int_\Omega (u_{x_1}+u_{x_2})v\,\mathrm dx
    -\int_\Omega fv\,\mathrm dx.
\end{equation}
The natural boundary condition $\Delta u=0$ on $\partial\Omega$ is incorporated through the weak formulation. To enforce the essential condition $u=0$ on $\partial\Omega$, we use the hard-constrained parameterization $ u_\theta=\mathcal N_\theta w$ with $ w(x)=4^{10}\prod_{i=1}^{10}x_i(1-x_i).$ The exact solution is chosen as $u^*(x)=\prod_{i=1}^{10}\sin(\pi x_i)$.

\begin{figure}
    \centering
    \includegraphics[width=0.75\linewidth]{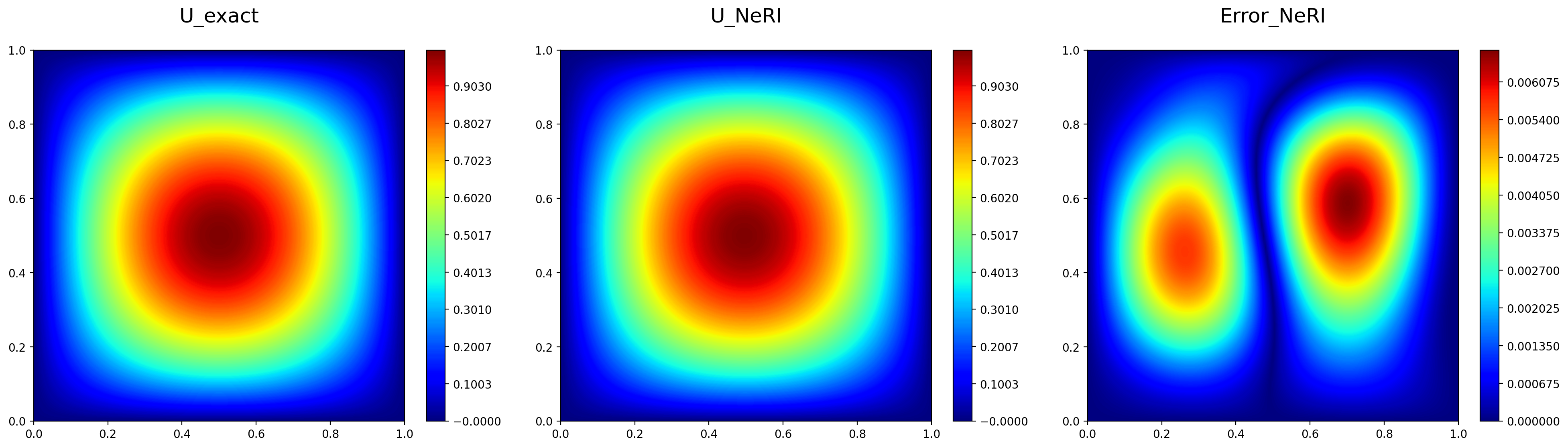}
    \caption{Numerical results for Example \ref{nume:biharmonic}.}
    \label{fig:biharmonic}
\end{figure}

The network is trained for $5{,}000$ epochs. The numerical solution on the two-dimensional slice
$(x_1,x_2,0.5,\ldots,0.5)$ is shown in Fig.~\ref{fig:biharmonic}. On this slice, NeRI achieves a relative $L^2$ error of $5.25\times10^{-3}$ and a maximum pointwise error of $6.63\times10^{-3}$.

\subsection{A nonsymmetric variational inequality problem}
\label{num:VI}

Let $\Omega=(-1,1)^2$. We consider the nonsymmetric variational inequality \eqref{eq:VI}, where
\begin{equation*}
    \langle \mathcal A u-f,v\rangle_{V^*,V}
    :=
    \int_\Omega
    \bigl(
        \nabla u\cdot\nabla v
        +u_{x_1}v
        +uv-fv
    \bigr)\,\mathrm dx
    -\int_{\partial\Omega}gv\,\mathrm ds,
\end{equation*}
and the constraint set is $K:=\{u\in H^1(\Omega): u\ge1 \text{ a.e. in }\Omega\}$. The exact solution is chosen as
\begin{equation*}
    u^*(x)=
    \begin{cases}
        (1-8x_1^2-8x_2^2)^2+1,
        & x_1^2+x_2^2\le\frac18,\\
        1,
        & \text{otherwise}.
    \end{cases}
\end{equation*}
The Neumann data $g$ are determined from $u^*$, and the source term is given by
\begin{equation*}
    f(x)
    =
    -\Delta u^*
    +u_{x_1}^*
    +u^*
    -
    \begin{cases}
        1, & x_1^2+x_2^2>\frac34,\\
        0, & \text{otherwise}.
    \end{cases}
\end{equation*}
The coincidence set is $\{u^*=1\}=\{x_1^2+x_2^2\ge 1/8\}$. In particular, the problem contains a biactive region where the obstacle is active while the multiplier $\mathcal A u^*-f$ vanishes. Moreover, the second derivatives of $u^*$ have jumps across the circle
$\{x_1^2+x_2^2=1/8\}$, resulting in limited global regularity.

\begin{figure}
    \centering
    \includegraphics[width=0.85\linewidth]{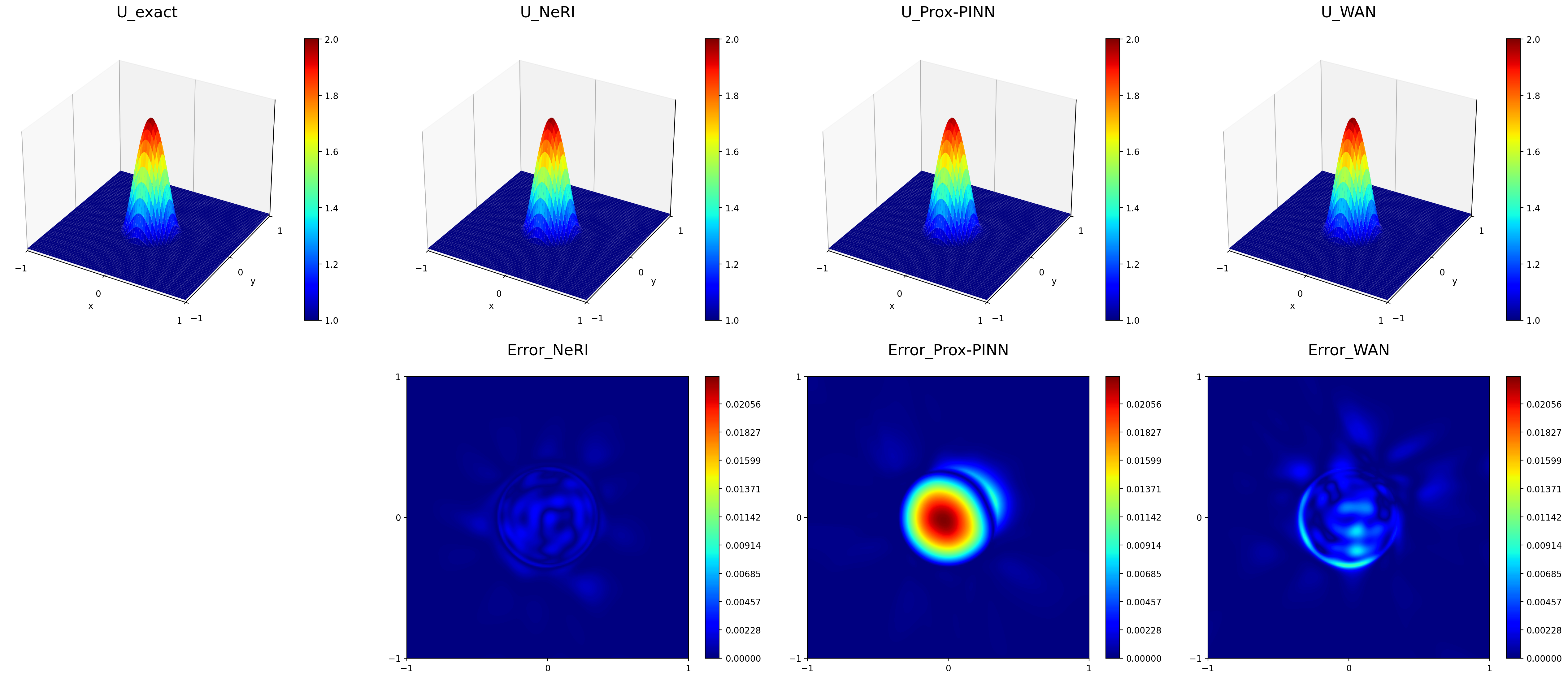}
    \caption{Numerical results for Example \ref{num:VI}.}
    \label{fig:Poisson}
\end{figure}

We compare NeRI with Prox-PINN \cite{gao2025prox} and the weak adversarial approach \cite{alphonse2024neural}. For all methods, we use the hard-constrained parameterization
$u_\theta=\mathcal N_\theta^2+1$
to enforce the obstacle constraint. For this problem, the Prox-PINN loss is
\begin{equation*}
    \mathcal L(\theta)
    =
    \frac{1}{\tau}
    \left\|
        \max\left\{
            u_\theta
            -\tau\left(
                \mathcal{A}u_\theta-f
            \right),
            1
        \right\}
        -u_\theta
    \right\|_{L^2(\Omega)}^2
    +\beta
    \left\|
        \frac{\partial u_\theta}{\partial n}-g
    \right\|_{L^2(\partial\Omega)}^2.
\end{equation*}
Compared with the standard Prox-PINN formulation, we include the factor $1/\tau$ in the first term, which improved the numerical accuracy in our parameter tuning. We use $\beta=10^2$ and $\tau=0.01$. For the weak adversarial method, we introduce two constraint-preserving networks $u_{\theta_u}=\mathcal{N}_{\theta_u}^2+1$ and $v_{\theta_v}=\mathcal{N}_{\theta_v}^2+1$, and use the objective
\begin{equation*}
\label{eq:vi-wan-regularized}
\mathcal{L}_{\mathrm{reg}}(\theta_u,\theta_v)
:=
\bigl\langle
Au_{\theta_u}-f,\,
u_{\theta_u}-v_{\theta_v}
\bigr\rangle_{V^*,V}
-
\gamma_{\mathrm{reg}}
\bigl\|u_{\theta_u}-v_{\theta_v}\bigr\|_{H^1(\Omega)}^{2},\qquad \gamma_{\mathrm{reg}}=5\times 10^{-4}.
\end{equation*}
The networks are trained for $20{,}000$ epochs. Fig.~\ref{fig:Poisson} shows the numerical approximations obtained by NeRI, Prox-PINN, and WAN together with the corresponding pointwise absolute errors. NeRI captures the localized nonsmooth behavior more accurately and yields the smallest error among the three methods. The relative $L^2$ errors for NeRI, Prox-PINN, and WAN are
$4.18\times10^{-4}$, $3.80\times10^{-3}$, and $1.09\times10^{-3}$, respectively.

We note that, although not presented here, for smooth problems with high-regularity $u^*$ (e.g., examples in \cite{alphonse2024neural}) the Prox-PINN can work quite well and tends to achieve higher accuracy and faster convergence than NeRI and WAN with a suitable choice of $\tau$ and $\beta$.

\section{Conclusion}
\label{sec:conclusion}

In this work, we introduced Neural Residual Iteration (NeRI), a weak-form neural-network framework for a class of elliptic operator equations and variational inequalities, including nonlinear and nonsymmetric problems. At each iteration, NeRI constructs a state-dependent linear functional whose derivative recovers the current weak residual. This construction yields a residual-driven parameter update without requiring a global energy functional or an additional test network.  We also extended NeRI to elliptic variational inequalities through penalized and constraint-preserving formulations and provided weak adversarial interpretations without explicit min--max training.

From a function-space perspective, NeRI induces, to first order, a neural-tangent-operator-preconditioned residual iteration. For strongly monotone Lipschitz operator equations, we established convergence of the fixed-kernel iteration to a reduced-space solution and developed a gradient-flow analysis in the Sobolev setting. For wide two-layer networks, we quantified the finite-time kernel drift in $\mathcal L(H^1(\Omega)^*,H^1(\Omega))$ and derived a direct $H^1(\Omega)$-error bound relative to the exact weak solution. This bound separates finite-width lazy-training, frozen-flow optimization, and reduced-space approximation errors, clarifying the role of the neural parameterization in both approximation and residual preconditioning.

Numerical experiments illustrate the effectiveness of NeRI for nonsymmetric and low-regularity problems, including high-dimensional equations and an obstacle-type variational inequality. The reported comparisons show advantages over the tested PINN-type and weak adversarial baselines, while retaining a single-network, single-step update structure.

Several directions remain for future work. One is to develop convergence theory for finite-width networks beyond the lazy-training regime and to extend the explicit Sobolev-space estimates to deep architectures. Another is to analyze the discrete-time neural iteration, accounting for the evolving tangent operator, finite step sizes, and quadrature or Monte Carlo approximation of the weak residual. The use of Adam in the numerical experiments also motivates extending the theory to different optimizers, including momentum and adaptive methods with history-dependent dynamics. Further questions concern convergence guarantees for constraint-preserving VI updates and extensions to more general monotone or nonsmooth operators, time-dependent problems, and coupled PDE systems.\\[1ex]

\noindent {\bf Acknowledgement:} This work is funded by the Deutsche Forschungsgemeinschaft (DFG, German Research Foundation) under Germany’s
Excellence Strategy – The Berlin Mathematics Research Center MATH+ (EXC-2046/1, project ID: 390685689). In particular, MH gratefully acknowledges support from the DFG Excellence Cluster MATH+ (EXC 2046) through the project AA5-2 ``Robust Multilevel Training of Artificial Neural Networks." This work is further supported by the DFG through the Research Grants Programme project ``Nonsmooth machine-learning-based computational models and methods for ultrasound computed tomography (MLUSCT)" (Project number: 581139020).

\bibliography{references}

\end{document}